\documentclass{article}

\usepackage{arxiv}

\usepackage[utf8]{inputenc} 
\usepackage[T1]{fontenc}    
\usepackage{hyperref}       
\usepackage{url}            
\usepackage{booktabs}       
\usepackage{amsfonts}       
\usepackage{nicefrac}       
\usepackage{microtype}      
\usepackage{lipsum}		
\usepackage{graphicx}
\usepackage[numbers]{natbib}
\usepackage{doi}
\usepackage{amsmath,amsfonts,amssymb}
\usepackage{mathtools}
\usepackage{graphicx}
\usepackage{setspace}
\usepackage{tocloft}
\usepackage{soul}
\usepackage{svg} 
\usepackage{siunitx} 
\usepackage{longtable} 
\usepackage{appendix} 
\usepackage[auth-lg]{authblk} 

\usepackage{xcolor}
\definecolor{c1}{rgb}{0.6,0,1} 

\title{A Modular Framework for Implicit 3D-0D Coupling in Cardiac Mechanics}

\author[a,c]{Aaron L. Brown}
\author[b,c,d]{Matteo Salvador}
\author[g,h]{Lei Shi}
\author[b,c,d]{Martin R. Pfaller}
\author[a]{Zinan Hu}
\author[f]{Kaitlin E. Harold}
\author[i]{Tzung Hsiai}
\author[g,*]{Vijay Vedula}
\author[a,b,c,d,e,*]{Alison L. Marsden}

\affil[a]{Department of Mechanical Engineering, Stanford University, Stanford, CA, USA}
\affil[b]{Institute for Computational and Mathematical Engineering, Stanford University, Stanford, CA, USA}
\affil[c]{Stanford Cardiovascular Institute, Stanford, CA, USA}
\affil[d]{Department of Pediatrics (Cardiology), Stanford University, Stanford, CA, USA}
\affil[e]{Department of Bioengineering, Stanford University, Stanford, CA, USA}
\affil[f]{Department of Computer Science, Stanford University, Stanford, CA, USA}
\affil[g]{Department of Mechanical Engineering, Columbia University, New York, NY, USA}
\affil[h]{Department of Mechanical Engineering, Kennesaw State University, Marietta, GA, USA}
\affil[i]{Department of Bioengineering, University of California Los Angeles, Los Angeles, CA, USA}
\affil[*]{Corresponding authors: vv2316@columbia.edu, amarsden@stanford.edu}

\hypersetup{
pdftitle={A Modular Framework for Implicit 3D-0D Coupling in Cardiac Mechanics},
pdfsubject={Cardiovascular modeling},
pdfauthor={Aaron L. Brown},
pdfkeywords={Cardiovascular modeling, cardiac mechanics, 3D-0D coupling, multi-domain modeling, Approximate Newton Method},
}

\begin{document}
\maketitle

\begin{abstract}
In numerical simulations of cardiac mechanics, coupling the heart to a model of the circulatory system is essential for capturing physiological cardiac behavior. A popular and efficient technique is to use an electrical circuit analogy, known as a lumped parameter network or zero-dimensional (0D) fluid model, to represent blood flow throughout the cardiovascular system. Due to the strong \textit{physical} interaction between the heart and the blood circulation, developing accurate and efficient \textit{numerical} coupling methods remains an active area of research. In this work, we present a modular framework for implicitly coupling three-dimensional (3D) finite element simulations of cardiac mechanics to 0D models of blood circulation. The framework is modular in that the circulation model can be modified independently of the 3D finite element solver, and vice versa. The numerical scheme builds upon a previous work that combines 3D blood flow models with 0D circulation models (3D fluid - 0D fluid). Here, we extend it to couple 3D cardiac tissue mechanics models with 0D circulation models (3D structure - 0D fluid), showing that both mathematical problems can be solved within a unified coupling scheme. The effectiveness, temporal convergence, and computational cost of the algorithm are assessed through multiple examples relevant to the cardiovascular modeling community. Importantly, in an idealized left ventricle example, we show that the coupled model yields physiological pressure-volume loops and naturally recapitulates the isovolumic contraction and relaxation phases of the cardiac cycle without any additional numerical techniques. Furthermore, we provide a new derivation of the scheme inspired by the Approximate Newton Method of Chan (1985), explaining how the proposed numerical scheme combines the stability of monolithic approaches with the modularity and flexibility of partitioned approaches. 
\end{abstract}

\keywords{Cardiovascular modeling \and cardiac mechanics \and 3D-0D coupling \and multi-domain modeling \and Approximate Newton Method}

\section{Introduction}\label{sect:intro}
Numerical simulations have long been used to investigate the cardiovascular system in both health and disease \cite{schwarz2023beyond}. These efforts have primarily applied computational fluid dynamics (CFD) to study blood flow in the heart and vasculature \cite{marsden2009evaluation, mittal2016computational, vedula2017method, zingaro2021geometric, mahmoudi2022guiding, nair2023hemodynamics}, computational solid dynamics (CSD) to simulate tissue mechanics in the heart and vasculature \cite{baillargeon2014living, strocchi2020simulating, karabelas2022accurate, fedele2023comprehensive, barnafi2023comparative}, and fluid-structure interaction (FSI) for coupled problems \cite{hirschhorn2020fluid, bucelli2022partitioned, bucelli2023mathematical, davey2023simulating}. Because accounting for the entire 3D circulatory system is typically infeasible due to limited imaging domains and a vast range of scales (micro- to macro-vessels), it is common to model parts of this system with a lumped parameter network (LPN), which can be thought of as a 0D model of blood flow. This treats blood flow in the circulatory system analogously to the flow of current in an electrical circuit \cite{quarteroni2016geometric} and allows one to quantify bulk quantities -- pressure and flow rate -- at various locations in the system, at a fraction of the cost of fully-resolved 3D CFD simulations. Representing some parts of the cardiovascular system -- for example, the heart or specific blood vessels -- with 3D structural and/or fluid models, while modeling the remainder using a 0D LPN constitutes a multi-domain approach \cite{vignon2006outflow, moghadam2013modular}. The 0D LPN acts as a boundary condition on the 3D model that recapitulates physiological effects not captured by other, simpler conditions (e.g., zero-pressure). Such a 3D-0D coupled problem is the focus of this work.

Developing accurate and efficient numerical methods for 3D-0D coupling remains an active area of research \cite{hirschvogel2017monolithic, augustin2021computationally, regazzoni2022cardiac, carichino2018energy}. Previous works have coupled 3D blood flow in large vessels to 0D models of the downstream vasculature \cite{vignon2006outflow, moghadam2013modular, quarteroni2016geometric}, while other groups have coupled 3D finite element models of the heart to 0D models of the systemic and pulmonary circulation \cite{jafari2019framework, shavik2018high, augustin2021computationally, regazzoni2022cardiac, zhang2023simulating}. While these two problems, 3D fluid - 0D fluid and 3D structure - 0D fluid, are related, to the best of our knowledge, none have previously treated them in a unified manner. 

Prior works have taken a variety of approaches to solving the coupled problem, which can be broadly categorized into monolithic \cite{vignon2006outflow, hirschvogel2017monolithic, augustin2021computationally} and partitioned \cite{kerckhoffs2007coupling, shavik2018high, regazzoni2022cardiac, carichino2018energy}. Monolithic schemes are robust and generally exhibit better convergence properties, but are not conducive to modularity. We define modular implementations as those in which the 3D and 0D equations are solved separately by independent codes optimized for their respective problems, and those codes exchange information as needed to couple the two sets of equations. Partitioned approaches are typically modular, but may suffer from numerical stability issues. For the problem of coupling 3D heart models to 0D circulation models, partitioned schemes suffer from a particular issue known as the balloon dilemma \cite{regazzoni2022cardiac}, which originates from the cardiac valves. During the two isovolumic phases of the cardiac cycle, both the inlet and outlet valves of the left ventricle (LV), for example, are closed, and the LV volume is nearly constant, while the LV pressure increases or decreases greatly due the contraction or relaxation of the heart muscle. In this situation, partitioned schemes that alternate between structure and fluid solvers typically fail because the structure solver is not aware of the constant-volume constraint. Previous works have avoided this by choosing monolithic approaches, or by using special iterative methods, time-staggered schemes, or additional stabilization terms.

In \cite{moghadam2013modular}, the authors developed a hybrid approach to the 3D fluid - 0D fluid coupling problem, incorporating the advantages of monolithic and partitioned approaches. This method was implemented in the open source multiphysics finite element solver svFSI (\url{https://github.com/SimVascular/svFSI}) \cite{zhu2022svfsi} and has been used extensively in blood flow simulations where the solution in the 3D domain is strongly influenced by the surrounding vascular system \cite{baumler2020fluid, lan2022virtual, gutierrez2019hemodynamic}. 

In this work, we describe a modular numerical scheme to implicitly couple 3D fluid and/or structural mechanics models to 0D LPNs of the cardiovascular system. The algorithm was originally described in \cite{moghadam2013modular} for only the 3D fluid - 0D fluid problem. Here, we extend it to solve the 3D structure - 0D fluid problem, showing that these two problems can be treated under a unified coupling framework. Applying this coupling to an idealized left ventricle model, we demonstrate the method produces a realistic pressure-volume loop and naturally captures the isovolumic cardiac phases and opening and closing of valves without additional numerical treatment to solve the balloon dilemma. We further derive the coupling scheme as a modification to the monolithic Newton approach, inspired by the Approximate Newton Method (ANM) of Chan \cite{chan1985approximate}, revealing a firm mathematical foundation. This connection to ANM also makes clear how the present coupling retains the robustness of a monolithic approach within a modular implementation like a partitioned approach. The modularity greatly improves usability, allowing the user to modify the 0D LPN independently of the 3D solver and vice versa. 

The paper is organized as follows. In Section \ref{sect:methods}, the proposed coupling framework is derived. In Section \ref{sect:results}, we leverage our numerical scheme in three different test cases, including an ellipsoidal LV coupled to an open-loop LPN, a spherical shell inflated through a limit point, and a pulmonary arterial model coupled to a closed-loop LPN. We also provide preliminary results on the convergence and computational cost of our method. In Section \ref{sect:discussion}, we consider our method in relation to recent works and discuss limitations and future directions. Finally, in Section \ref{sect:conclusion}, we summarize our findings with respect to the proposed scheme.

\section{Methods} \label{sect:methods}

In this section, we derive the proposed coupling framework. The resulting equations are the same as those given in \cite{moghadam2013modular, esmaily2013new}, but a generalized mathematical derivation, inspired by ANM \cite{chan1985approximate} and applicable to both 3D fluid and 3D structure problems, is provided here.  First, we state the governing equations for the 3D and 0D systems, then we describe how the two systems are mathematically coupled, and finally, we explain how to solve the coupled problem in a modular manner. In the following, the minor differences in the equations when considering a 3D fluid versus a 3D structure are highlighted.

\subsection{3D mechanical model: fluid or structure}
\begin{figure} [h]
    \centering
    \includegraphics[width = 0.9 \textwidth]{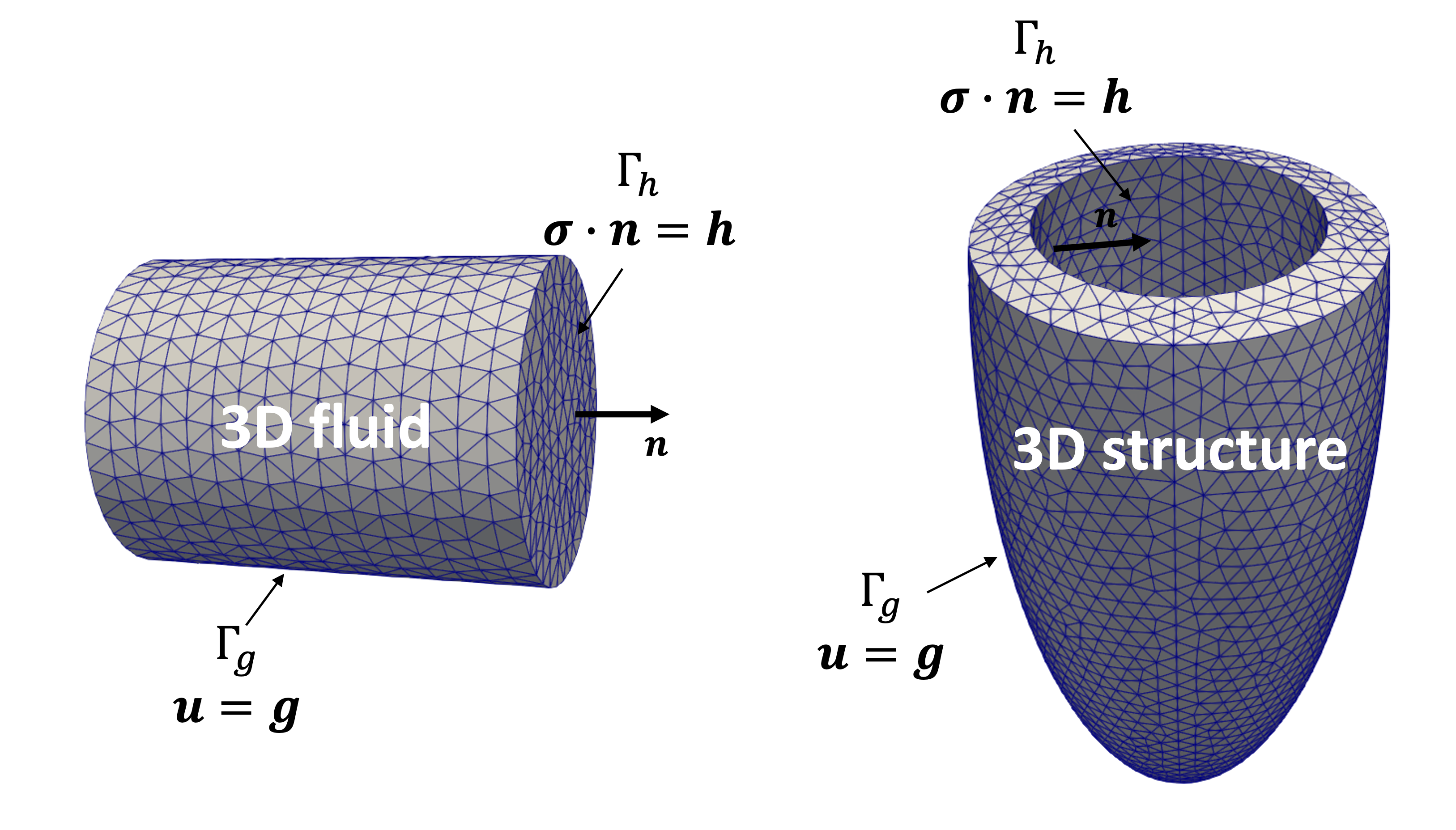}
    \caption{Left: An idealized geometry of a section of a blood vessel is given as an example 3D fluid domain. Right: An idealized geometry of the LV of the heart is given as an example 3D structure domain. The dynamics in each are described by standard governing partial differential equations (PDEs), and on both, we may define Dirichlet or Neumann boundary conditions, or a combination of both. In this work, the PDEs are spatially discretized using the finite element method.
    }
    \label{fig:3D_model}
\end{figure}
We first state the governing equations and numerical formulation for an incompressible and Newtonian fluid on a fixed 3D domain \cite{moghadam2013modular}, which models blood flow in large blood vessels (Fig. \ref{fig:3D_model} left). Specifically, the Navier-Stokes equations, consisting of the momentum and continuity equations, as well as the Newtonian constitutive model, read
\begin{equation} \label{eq:NS_momentum}
    \rho \frac{\partial \mathbf{u}}{\partial t} + \rho \mathbf{u} \cdot \nabla \mathbf{u} - \nabla \cdot \boldsymbol{\sigma} - \mathbf{f} = \mathbf{0} ,
\end{equation}
\begin{equation} \label{eq:NS_continuity}
    \nabla \cdot \mathbf{u} = 0,
\end{equation}
\begin{equation}
    \boldsymbol{\sigma} = -p \mathbf{I} + \mu (\nabla \mathbf{u} + \nabla \mathbf{u}^T),
\end{equation}
with boundary conditions
\begin{equation} \label{eq:DirichletBC}
    \mathbf{u} = \mathbf{g}, \hspace{5pt} \mathbf{x} \in \Gamma_g,
\end{equation}
\begin{equation} \label{eq:NeumannBC}
    \boldsymbol{\sigma} \cdot \mathbf{n} = \mathbf{h}, \hspace{5pt} \mathbf{x} \in \Gamma_h,
\end{equation}
and initial conditions
\begin{equation}
    \mathbf{u}(t=0) = \mathbf{u}_0,
\end{equation}
\begin{equation}
    p(t=0) = p_0,
\end{equation}
with position vector $\mathbf{x}$, time $t$, density $\rho$, velocity $\mathbf{u}$,  Cauchy stress tensor $\boldsymbol{\sigma}$, pressure $p$, dynamic viscosity $\mu$, body force $\mathbf{f}$, which we assume to be zero, and surface normal vector $\mathbf{n}$. Eq. \eqref{eq:DirichletBC} is a Dirichlet boundary condition with prescribed velocity $\mathbf{g}$ on $\Gamma_g$. Similarly, Eq. \eqref{eq:NeumannBC} is a Neumann boundary condition with prescribed traction $\mathbf{h}$ on $\Gamma_h$. We may write these equations in abstract form as

\begin{equation}
    \begin{cases}
        \mathcal{P}^{3D, fluid}(\mathbf{u},p, \mathbf{x}, t) = \mathbf{0},
        \\
        \text{Boundary conditions},
        \\
        \text{Initial conditions}.
    \end{cases}
\end{equation}

Following \cite{moghadam2013modular}, these equations are discretized in space using a stabilized (variational multiscale) finite element formulation and in time using the generalized-$\alpha$ method. This yields the nonlinear residual equation at timestep $n+1$
\begin{equation}
    \mathbf{R}^{3D,fluid}(\dot{\mathbf{U}}_{n+1}, \boldsymbol{\Pi}_{n+1}) = \mathbf{0},
\end{equation}
to be solved for $\dot{\mathbf{U}}_{n+1}$ and $\boldsymbol{\Pi}_{n+1}$, the vectors of nodal accelerations and nodal pressures at the next timestep $n+1$, respectively. The residual is also a function of $\dot{\mathbf{U}}_{n}$ and  $\boldsymbol{\Pi}_{n}$, but these are assumed to be known, and thus we do not explicity write the functional dependence on them. This equation is solved using Newton's method, which in turn requires solving the following linear system at each Newton iteration $k$
\begin{equation} \label{eq:fluid_linear_system}
    \begin{bmatrix}
    \tilde{\mathbf{K}} & \mathbf{G} \\
    \mathbf{D} & \mathbf{L}
    \end{bmatrix}^{(k)}_{n+1}
    \begin{bmatrix}
    \Delta \dot{\mathbf{U}}^{(k)}_{n+1} \\
    \Delta \boldsymbol{\Pi}^{(k)}_{n+1}
    \end{bmatrix}
    =
    -
    \begin{bmatrix}
    \mathbf{R}^{3D,fluid}_m \\
    \mathbf{R}^{3D,fluid}_c
    \end{bmatrix}^{(k)}_{n+1}
    .
\end{equation}
$\mathbf{R}^{3D,fluid}_m$ is the residual associated with momentum balance Eq. \eqref{eq:NS_momentum}, while $\mathbf{R}^{3D,fluid}_c$ is the residual associated with mass continuity Eq. \eqref{eq:NS_continuity}. $\tilde{\mathbf{K}}, \mathbf{G},
    \mathbf{D}, \mathbf{L}$ are blocks of the tangent or stiffness matrix, and $\Delta \dot{\mathbf{U}}^{(k)}_{n+1}$ and $
    \Delta \boldsymbol{\Pi}^{(k)}_{n+1}$ are the Newton increments in nodal accelerations and pressures, respectively. The notation 
$\begin{bmatrix}
    \cdot
\end{bmatrix}^{(k)}_{n+1}$
indicates that terms inside the brackets are evaluated at $\dot{\mathbf{U}}_{n+1}^{(k)}$ and $\boldsymbol{\Pi}_{n+1}^{(k)}$. The solutions are updated each Newton iteration until convergence according to
\begin{equation}
    \dot{\mathbf{U}}_{n+1}^{(k+1)} =  \dot{\mathbf{U}}_{n+1}^{(k)} + \Delta \dot{\mathbf{U}}^{(k)}_{n+1},
\end{equation}
\begin{equation}
    \boldsymbol{\Pi}_{n+1}^{(k+1)} =  \boldsymbol{\Pi}_{n+1}^{(k)} + \Delta \boldsymbol{\Pi}_{n+1}^{(k)}.
\end{equation}

In cardiovascular biomechanics modeling, we are not only interested in blood flow, but also in the dynamics of the tissues surrounding the blood, notably the heart (Fig. \ref{fig:3D_model} right). The deformation of these tissues is governed by the equations of finite deformation elastodynamics. Specifically, we may state the Cauchy momentum equation in Lagrangian form
\begin{equation}
   \rho \frac{D \mathbf{u}}{D t} - \nabla \cdot \boldsymbol{\sigma} - \mathbf{f} = \mathbf{0},
\end{equation}
where $\frac{D}{Dt}$ denotes the material derivative. As with the fluid equations, $\mathbf{u}$ is the velocity, $\boldsymbol{\sigma}$ is the Cauchy stress tensor, and $\mathbf{f}$ is a body force, which we assume to be zero in this work. In our finite element solver, the structural problem is solved in the reference configuration, in which the relevant stress measure is the second Piola-Kirchhoff stress
\begin{equation}
    \mathbf{S} = J \mathbf{F}^{-1} \boldsymbol{\sigma} \mathbf{F}^{-T},
\end{equation}
where $\mathbf{F}$ is the deformation gradient tensor and $J = \det{\mathbf{F}}$ is the Jacobian. For a hyperelastic material described by a strain energy density function $\psi(\mathbf{F})$, we have
\begin{equation}
    \mathbf{S}
    =\frac{\partial \psi}{\partial \mathbf{E}},
\end{equation}
where $\mathbf{E} = \frac{1}{2} (\mathbf{C} - \mathbf{I})$ is the Green-Lagrange strain tensor and $\mathbf{C} = \mathbf{F}^T \mathbf{F}$ is the right Cauchy-Green tensor. 
These equations are augmented with boundary and initial conditions so that we may write the structural dynamics problem in abstract form
\begin{equation}
    \begin{cases}
        \mathcal{P}^{3D, struct}(\mathbf{u}, \mathbf{x}, t) = \mathbf{0},
        \\
        \text{Boundary conditions},
        \\
        \text{Initial conditions}.
    \end{cases}
\end{equation}

These equations are solved using similar methods to those for fluid flow (i.e., finite element method and generalized-$\alpha$ method) \cite{hirschvogel2017monolithic}, leading to an analogous nonlinear system of equation to be solved at each timestep
\begin{equation}
    \mathbf{R}^{3D, struct}(\dot{\mathbf{U}}_{n+1}) = \mathbf{0}
    .
\end{equation}
Usually, nodal displacements are chosen as the structural unknowns after time discretization, but in our implementation, we instead choose nodal accelerations. This is an arbitrary choice that allows the structure problem to be treated similarly to the fluid problem, but it is not necessary for the present coupling framework. 
The system is solved using Newton's method, 
\begin{equation}
    \begin{bmatrix}
    \mathbf{K}
    \end{bmatrix}
    ^{(k)}_{n+1}
    \begin{bmatrix}
    \Delta \dot{\mathbf{U}}^{(k)}_{n+1} 
    \end{bmatrix}
    = -
    \begin{bmatrix}
    \mathbf{R}^{3D, struct}
    \end{bmatrix}
    ^{(k)}_{n+1},
\end{equation}
where $k$ again indicates the Newton iteration.

In the remainder of the paper, we consider the general 3D problem
\begin{equation}
    \begin{cases}
        \mathcal{P}^{3D}(\boldsymbol{\phi}, \mathbf{x}, t) = \mathbf{0},
        \\
        \text{Boundary conditions},
        \\
        \text{Initial conditions},
    \end{cases}
\end{equation}
where $\mathcal{P}^{3D}$ may represent either the 3D fluid or 3D structure PDE. $\boldsymbol{\phi}$ are the 3D variables, which include velocity or velocity and pressure, depending on the physics. After space-time discretization, we obtain the general 3D residual equation as 
\begin{equation} \label{eq:3D_res}
    \mathbf{R}^{3D}(\boldsymbol{\Phi}_{n+1}) = \mathbf{0}.
\end{equation}
$\boldsymbol{\Phi}_{n+1}$ is the state vector of the 3D system, where, for the fluid, $\boldsymbol{\Phi}_{n+1} = \begin{bmatrix}
    \dot{\mathbf{U}}_{n+1}
    \\
    \mathbf{\Pi}_{n+1}
\end{bmatrix}$,
while for the structure, $\boldsymbol{\Phi}_{n+1} = \dot{\mathbf{U}}_{n+1}$. Newton's method to solve this nonlinear system gives the general 3D linear system
\begin{equation} \label{eq:3D_linear_system}
    \begin{bmatrix}
    \displaystyle
    \frac{\partial \mathbf{R}^{3D}}{\partial \boldsymbol{\Phi}_{n+1}}
    \end{bmatrix}_{n+1}^{(k)}
    \begin{bmatrix}
    \Delta \boldsymbol{\Phi}_{n+1}^{(k)}
    \end{bmatrix}
    =
    -
    \begin{bmatrix}
    \mathbf{R}^{3D}
    \end{bmatrix}_{n+1}^{(k)}.
\end{equation}

\subsection{0D circulation model}
\begin{figure} [h]
    \centering
    \includegraphics[width = 0.3\linewidth]{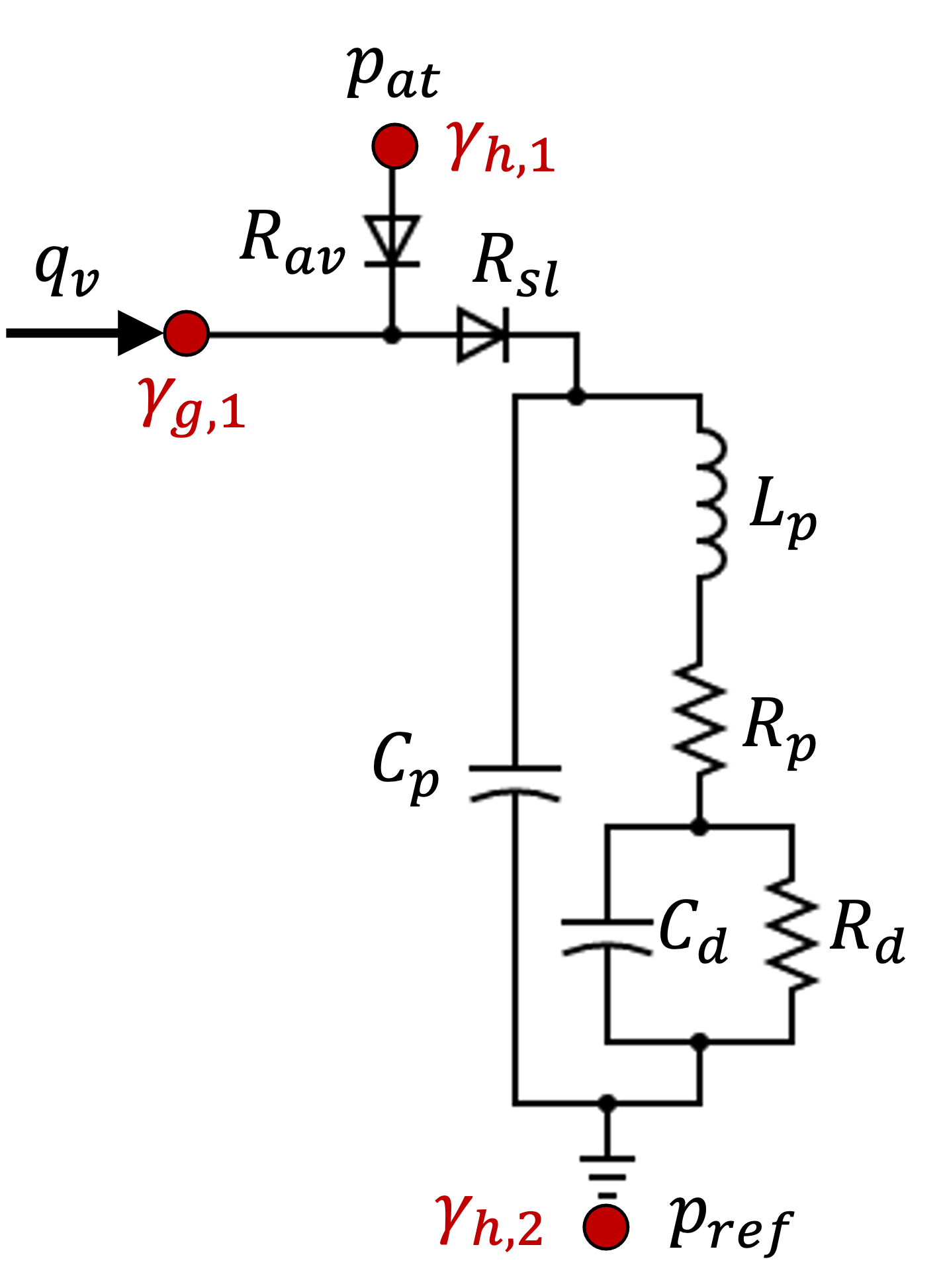}
    \caption{An example of a 0D circulation model or LPN, which treats blood flow through the body like the flow of current through an electrical circuit. The LPN is forced by prescribed flow rates $\mathbf{q}$ at the Dirichlet boundary nodes $\gamma_g = \{ \gamma_{g,1}\}$ and prescribed pressures $\mathbf{p}$ at the Neumann boundary nodes $\gamma_h = \{\gamma_{h,1}, \gamma_{h,2} \}$, both of which are denoted by red circles.}
    \label{fig:0D_model}
\end{figure}
Blood flow throughout the circulatory system is modeled using an LPN \cite{hirschvogel2017monolithic, regazzoni2022cardiac, augustin2021computationally}, also called a 0D fluid model. Fig. \ref{fig:0D_model} shows an example LPN that will be used later in Section \ref{sect:LV_sim}. LPN models are typically combinations of resistors $R$, which model the viscous resistance of vessels to blood flow, inductors $L$, which account for the inertia of blood, and capacitors $C$, which model the compliance of blood vessels. In addition, diodes are used to represent the heart valves. Many LPN models of the vasculature have been used in the literature, ranging from simple resistance or resistance-capacitance models of arteries to extensive closed-loop networks representing the entire circulatory system \cite{moghadam2013modular, regazzoni2022cardiac, kung2013predictive}. 

An LPN can be analyzed using Kirchhoff's first law for an electrical circuit, which leads to a general representation as a system of differential-algebraic equations (DAEs) \cite{peiro2009reduced},
\begin{equation} \label{eq:0D_ode}
    \frac{d \mathbf{y}}{dt}= \mathbf{f}(\mathbf{y}, \mathbf{z}, t),
\end{equation}
\begin{equation} \label{eq:0D_algebraic}
    \mathbf{g}(\mathbf{y}, \mathbf{z}, t) = \mathbf{0},
\end{equation}
with initial conditions
\begin{equation} \label{eq:0D_init_cond}
    \mathbf{y}(t = 0) = \mathbf{y}_0,
\end{equation}
where $\mathbf{y}$ contains differential variables determined by the differential equations Eq. \eqref{eq:0D_ode},  $\mathbf{z}$ contains algebraic variables determined by the algebraic equations Eq. \eqref{eq:0D_algebraic}, and $t$ is time. Both $\mathbf{y}$ and $\mathbf{z}$ contain pressures and flow rates at nodes and branches, respectively, and may also contain other variables, such as the cross-sectional area of a vessel or the state of a valve. $\mathbf{f}$ and $\mathbf{g}$ are a potentially nonlinear functions of $\mathbf{y}$, $\mathbf{z}$, and $t$. 

As with the 3D model, we may define boundary conditions for the 0D-LPN model. Let $\gamma_g$ denote the set of Dirichlet boundary nodes on which we prescribe flows $\mathbf{q}$. Let $\gamma_h$ denote the set of Neumann boundary nodes on which we prescribe pressures $\mathbf{p}$. Fig. \ref{fig:0D_model} gives an example LPN with these boundary nodes shown in red. The boundary flow rates and pressures are not boundary conditions in a strict sense because there is no notion of length in a 0D fluid model. Instead, they act as forcing terms directly in the 0D equations \cite{quarteroni2016geometric} 
\begin{equation} \label{eq:0D_ode_with_q}
    \frac{d \mathbf{y}}{dt} = \mathbf{f}(\mathbf{y}, \mathbf{z}, t, \mathbf{q}, \mathbf{p}),
\end{equation}
\begin{equation} \label{eq:0D_algebraic_with_q}
    \mathbf{g}(\mathbf{y}, \mathbf{z}, t, \mathbf{q}, \mathbf{p}) = \mathbf{0}.
\end{equation}
$\mathbf{y}$ and $\mathbf{z}$ can often be considered together, so we define the combined 0D state vector $\mathbf{w} = [\mathbf{y}, \mathbf{z}]^T$.

Analogous to the 3D system, the 0D equations are written in the following abstract form
\begin{equation}
    \begin{cases}
        \mathcal{P}^{0D}(\mathbf{w}, t, \mathbf{q}, \mathbf{p}) = \mathbf{0},
        \\
        \text{Initial conditions}.
    \end{cases}
\end{equation}

In this work, we integrate the 0D system with a 4th-order Runge-Kutta (RK4) scheme (Appendix \ref{sect:RK4_scheme}). We first apply RK4 to integrate the differential variables $\mathbf{y}$ using Eq. \eqref{eq:0D_ode_with_q} from timestep $n$ to $n+1$, then determine the algebraic variables $\mathbf{z}$ with Eq. \eqref{eq:0D_algebraic_with_q} using the updated differential variables. 
This yields a system of algebraic equations to be solved at timestep $n+1$
\begin{equation}
    \mathbf{R}^{0D}(\mathbf{w}_{n+1}) = \mathbf{0}.
\end{equation}
If the time-stepping scheme is explicit, as in this work, this system can be solved directly (i.e., in one iteration). If the scheme is implicit and the DAE system is nonlinear in $\mathbf{w}$, this system can be solved by Newton's method.

\subsection{The coupled problem} \label{sect:coupled_problem}
\begin{figure}[h]
    \centering
    \includegraphics[width = 0.9 \textwidth]{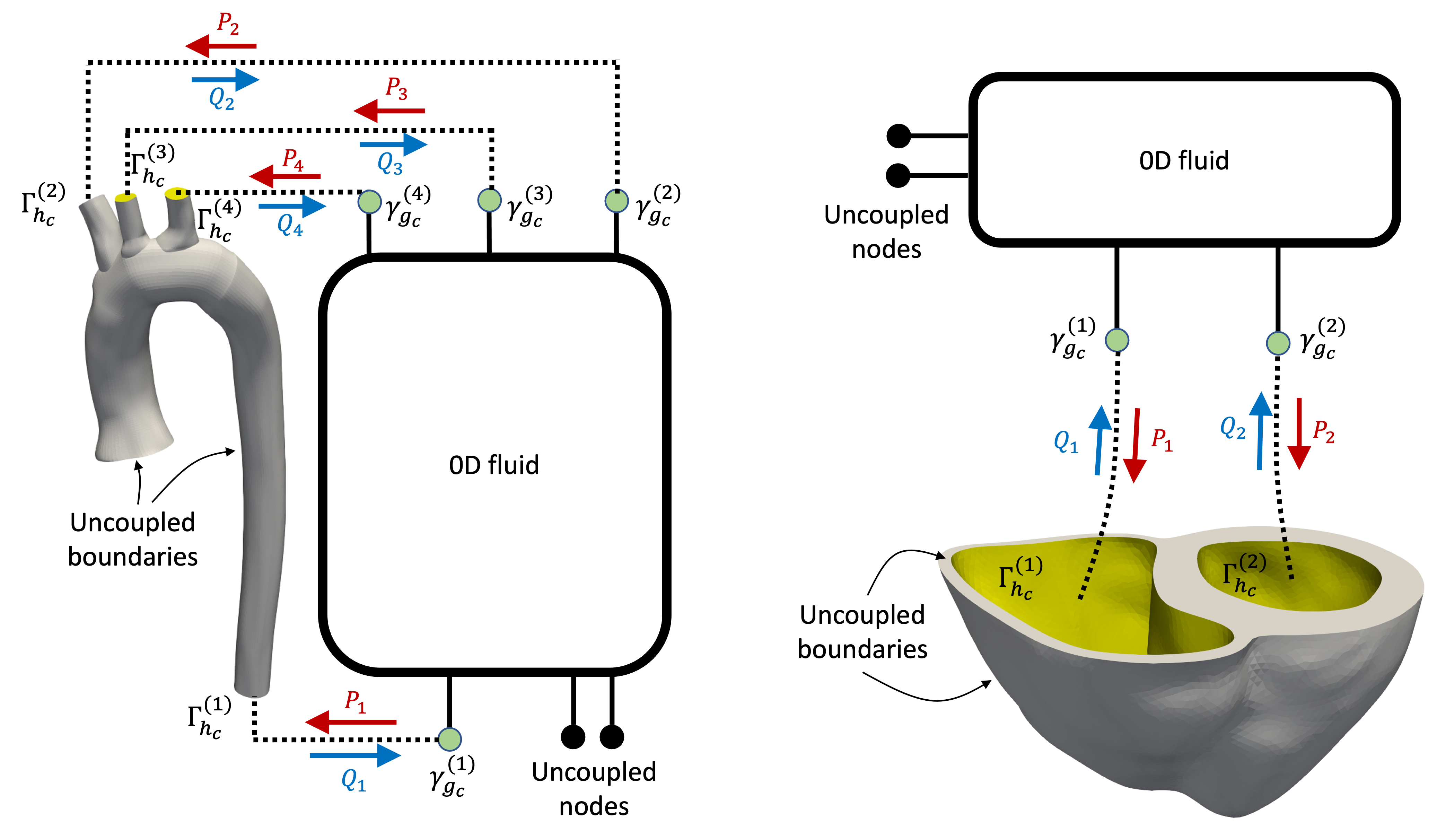}
    \caption{Left: Coupling between a 3D fluid and 0D fluid. The 3D fluid is an aorta model taken from the Vascular Model Repository (\url{https://www.vascularmodel.com/}). Right: Coupling between a 3D structure and 0D fluid. The 3D structure is a biventricular model obtained from patient MRI data. Both coupling problems are treated identically. Coupled Neumann surfaces $\Gamma_{h_c}^{(i)}$ on the 3D models are highlighted in yellow. Coupled Dirichlet nodes $\gamma_{g_c}^{(i)}$ on the 0D model are shown in green. Along each $\Gamma_{h_c}^{(i)}$ - $\gamma_{g_c}^{(i)}$ connection is an associated exchange of flow rate $Q_i$ and pressure $P_i$. Note that in addition to coupled boundaries, the 3D model will generally have additional uncoupled boundaries on which one may prescribe uncoupled Dirichlet and/or Neumann boundary conditions, and likewise the 0D model will generally have additional uncoupled nodes on which one may prescribe uncoupled Dirichlet and/or Neumann boundary forcings.
    }
    \label{fig:3D_0D_coupling}
\end{figure}
So far, we have individually discussed the numerical treatment of the 3D mechanical model and the 0D circulation model. In this section, we describe how these two models are mathematically coupled. 

The 3D Neumann boundary is split into coupled and uncoupled parts, $\Gamma_h = \Gamma_{h_c} \cup \Gamma_{h_u}$. 
In general, there may be multiple distinct coupled Neumann boundaries, so that $\Gamma_{h_c} = \Gamma_{h_c}^{(1)} \cup \Gamma_{h_c}^{(2)} \cup \ldots \cup \Gamma_{h_c}^{(n^{cBC})}$, where $n^{cBC}$ is the number of coupled boundaries. In Fig. \ref{fig:3D_0D_coupling}, these are the outflow boundaries of the aorta model and the endocardial surfaces of the biventricular model (i.e., the inner surfaces of the heart muscle in contact with the blood).
Similarly, the set of 0D Dirichlet nodes is split into coupled and uncoupled parts, $\gamma_g = \gamma_{g_c} \cup \gamma_{g_u}$. The set of coupled Dirichlet nodes may be written $\gamma_{g_c} = \gamma_{g_c}^{(1)} \cup \gamma_{g_c}^{(2)} \cup \ldots \cup
\gamma_{g_c}^{(n^{cBC})}$. Each 0D coupled Dirichlet node $\gamma_{g_c}^{(i)}$ is associated with a 3D coupled Neumann boundary $\Gamma_{h_c}^{(i)}$, $i \in \{1, \ldots, n^{cBC}\}$, as shown in Fig. \ref{fig:3D_0D_coupling}.

With these definitions in place, we may state the mathematical coupling between the 3D and 0D domains. For the 3D, we impose a spatially uniform pressure $P_i$ on the coupled 3D boundary $\Gamma_{h_c}^{(i)}$, where the value of $P_i$ is taken as the pressure at the corresponding coupled node $\gamma_{g_c}^{(i)}$ of the 0D model. Stated mathematically,

\begin{equation} \label{eq:P_relation}
    P_i = (\mathbb{P}\mathbf{w})_i \hspace{10pt} \text{on} \hspace{5pt} \Gamma_{h_c}^{(i)},
\end{equation}
where $\mathbb{P}$ is a matrix that selects the appropriate components from $\mathbf{w}$. For example, if $\mathbf{w}$ has 5 components, but only the second and fourth components represent pressures at coupled nodes 1 and 2, then
\begin{equation}
    \mathbb{P}
    =
    \begin{bmatrix}
    0 & 1 & 0 & 0 & 0 \\
    0 & 0 & 0 & 1 & 0
    \end{bmatrix}
    .
\end{equation}
Note that the spatially uniform pressure assumption is made in several other coupling approaches \cite{hirschvogel2017monolithic,regazzoni2022cardiac, augustin2021computationally}. 

Analogously, for the 0D, we impose a flow rate $Q_i$ at the coupled 0D node $\gamma_{g_c}^{(i)}$, where the value of $Q_i$ is taken as the velocity flux through the corresponding coupled 3D boundary $\Gamma_{h_c}^{(i)}$ of the 3D model. Stated mathematically,
\begin{equation} \label{eq:Q_relation}
     Q_i = \int_{\Gamma_{h_c}^{(i)}} \mathbf{u} \cdot \mathbf{n} d\Gamma,
\end{equation}
where $\mathbf{n}$ is the outward surface normal vector. This definition of flow rate should be clear in the context of a 3D fluid. For a 3D structure, however, in order to define a flow rate $Q_i$, we must restrict our attention to 3D structures that enclose some volume of blood $V_i$, most commonly a chamber of the heart. Then, we may define flow rate as $Q_i = -\frac{dV_i}{dt}$. In Appendix \ref{sect:flow_rate}, using the Reynolds Transport Theorem, it is shown that Eq. \eqref{eq:Q_relation} is equally valid for such a 3D structure, provided $\Gamma_{h_c}^{(i)}$ is a surface that closes the volume of interest $V_i$ (see Section \ref{sect:capping} if surface is not closed) and the integral is taken over $\Gamma_{h_c}^{(i)}$ in the deformed configuration. It is this observation that permits uniform treatment of 3D fluid and structure. 

On a practical note, in the finite element setting, the flow rate integral is computed from the 3D degrees of freedom as follows.
\begin{equation} \label{eq:Q3D}
    Q_i = \int_{\Gamma_{h_c}^{(i)}} \mathbf{u} \cdot \mathbf{n} d\Gamma = \sum_A \int _{\Gamma_{h_c}^{(i)}} N_A (\mathbf{U})_A \cdot \mathbf{n} d\Gamma,
\end{equation}
where $(\mathbf{U})_A$ is the velocity of node $A$ of the finite element model and $N_A (\mathbf{x})$ is the associated shape function in the finite element formulation. Note also that in the time-discrete setting, the nodal velocities are obtained from the nodal accelerations by the generalized-$\alpha$ method expression
\begin{equation} \label{eq:genalpha_U}
    \mathbf{U}_{n+1} = \mathbf{U}_n + \Delta t \dot{\mathbf{U}}_{n} + \gamma \Delta t (\dot{\mathbf{U}}_{n+1} - \dot{\mathbf{U}}_n),
\end{equation}
where $\Delta t$ is the timestep size and $\gamma$ is a parameter of the generalized-$\alpha$ method, not to be confused with the 0D Dirichlet and Neumann boundary node sets $\gamma_g$ and $\gamma_h$.
Depending on the physical formulation (structure or fluid), $\dot{\mathbf{U}}$ is either precisely $\boldsymbol{\Phi}$ or a component of $\boldsymbol{\Phi}$.

The coupled problem in the time and space continuous domain may be summarized in the following abstract manner
\begin{equation} \label{eq:abstract_coupled_problem}
    \begin{cases}
    \mathcal{P}^{0D}(\mathbf{w}, t, [\mathbf{q}_{u}, \mathbf{Q}], \mathbf{p}) = \mathbf{0},
    \\
    \text{Initial conditions},
    \\
    \\
    \mathcal{P}^{3D}(\boldsymbol{\phi}, \mathbf{x}, t) = \mathbf{0},
    \\
    \text{Initial conditions},
    \\
    \text{Uncoupled boundary conditions,}
    \\
    \boldsymbol{\sigma} \cdot \mathbf{n} = -P_i \mathbf{n}, \text{ on } \Gamma_{h_c}^{(i)}, i \in \{1, \ldots, n^{cBC}\},
    \\
    \\
   P_i(t) = (\mathbb{P}\mathbf{w}(t))_i,
    \\
    Q_i(t) = \int_{\Gamma_{h_c}^{(i)}(t)} \mathbf{u}(t) \cdot \mathbf{n}(t) d\Gamma,
    
    \end{cases}
\end{equation}
where the 0D Dirichlet forcing term $\mathbf{q}$ is split into an uncoupled component $\mathbf{q}_u$, which are prescribed flow rates on the 0D model, and a coupled component $\mathbf{Q}$, which is obtained from the 3D model (i.e., $\mathbf{q} = [\mathbf{q}_u, \mathbf{Q}]$). Similarly, the 3D boundary conditions are split into an uncoupled component (Dirichlet or Neumann), which is prescribed, and coupled pressure boundary conditions with magnitude $P_i$, which are obtained from the 0D model. The expressions for $P_i$ and $Q_i$, given in Eqs. \eqref{eq:P_relation} and \eqref{eq:Q_relation} and restated here, provide the coupling conditions between 0D and 3D. While the values of the uncoupled boundary conditions for both 3D and 0D are generally prescribed, the values of the coupled pressure and flow boundary conditions, $P_i$ and $Q_i$, are unknown and must be determined as part of the solution to the coupled problem.

After applying RK4 time discretization to the 0D system and applying generalized-$\alpha$ time discretization and finite element spatial discretization to the 3D system, the problem reduces to solving the following coupled nonlinear equations at each timestep $n+1$,
\begin{equation} \label{eq:coupled_residuals}
    \begin{cases}
    \mathbf{R}^{0D}(\mathbf{w}_{n+1},\mathbf{Q}_{n+1}(\boldsymbol{\Phi}_{n+1})) = \mathbf{0}
    ,
    \\
    \mathbf{R}^{3D}(\boldsymbol{\Phi}_{n+1}, \mathbf{P}_{n+1}(\mathbf{w}_{n+1})) = \mathbf{0}
    ,
    \end{cases}
\end{equation}
where, due to the coupling between 3D and 0D, the 0D residual $\mathbf{R}^{0D}$ is a function of the 3D state vector $\boldsymbol{\Phi}_{n+1}$ through the interface flow rates $\mathbf{Q}_{n+1}$, and similarly the 3D residual $\mathbf{R}^{0D}$ is a function of the 0D state vector $\mathbf{w}_{n+1}$ through the interface pressures $\mathbf{P}_{n+1}$.

\subsection{Solving the coupled problem}
The method to solve the system Eq. \eqref{eq:coupled_residuals} is inspired by ANM \cite{chan1985approximate}, with modifications suggested in \cite{artlich1995newton} and in \cite{matthies2003partitioned}. The coupled problem is solved in a modular manner, which allows us to take advantage of codes that already exist to efficiently solve the 3D and 0D problems. Moreover, the 0D model may be modified without changing the 3D solver, and vice versa. Using ANM as a foundation, we reproduce the equations of the original coupling framework described in \cite{moghadam2013modular, esmaily2013new}. In addition, we show the equations apply not only to 3D fluid - 0D fluid coupling, but also to 3D structure - 0D fluid coupling.

We begin by applying Newton's method to solve Eq. \eqref{eq:coupled_residuals} in a monolithic manner, which yields a linear system to be solved at each Newton iteration, $k$,
\begin{equation} \label{eq:monolithic_system}
    \begin{bmatrix}
    \displaystyle \frac{\partial \mathbf{R}^{0D}}{\partial \mathbf{w}_{n+1}} & 
    \displaystyle \frac{\partial \mathbf{R}^{0D}}{\partial \boldsymbol{\Phi}_{n+1}}
    \\
    \displaystyle \frac{\partial \mathbf{R}^{3D}}{\partial \mathbf{w}_{n+1}} & 
    \displaystyle \frac{\partial \mathbf{R}^{3D}}{\partial \boldsymbol{\Phi}_{n+1}}
    \end{bmatrix}^{(k)}_{n+1}
    \\
    \begin{bmatrix}
    \Delta \mathbf{w}^{(k)}_{n+1}
    \\
    \Delta \boldsymbol{\Phi}^{(k)}_{n+1}
    \end{bmatrix}
    =
    -
    \begin{bmatrix}
    \mathbf{R}^{0D}
    \\
    \mathbf{R}^{3D}
    \end{bmatrix}^{(k)}_{n+1},
\end{equation}
with the update
\begin{equation}
    \boldsymbol{\Phi}_{n+1}^{(k+1)} = \boldsymbol{\Phi}_{n+1}^{(k)} + \Delta \boldsymbol{\Phi}_{n+1}^{(k)}  \hspace{15pt}\text{and}\hspace{15pt}
    \mathbf{w}_{n+1}^{(k+1)} = \mathbf{w}_{n+1}^{(k)} + \Delta \mathbf{w}_{n+1}^{(k)}.
\end{equation}
As before, the notation 
$\begin{bmatrix}
    \cdot
\end{bmatrix}^{(k)}_{n+1}$
indicates that terms inside the brackets are evaluated at $\boldsymbol{\Phi}_{n+1}^{(k)}$ and $\mathbf{w}_{n+1}^{(k)}$.
Performing Schur Complement Reduction \cite{benzi2005numerical} (also known as Block Gauss Elimination or Static Condensation) yields the equivalent system
\begin{multline} \label{eq:ANM_block_gauss}
    \begin{bmatrix}
    \displaystyle \frac{\partial \mathbf{R}^{0D}}{\partial \mathbf{w}_{n+1}} 
    & 
    \displaystyle \frac{\partial \mathbf{R}^{0D}}{\partial \boldsymbol{\Phi}_{n+1}}
    \\
    \mathbf{0} 
    & 
    \displaystyle \frac{\partial \mathbf{R}^{3D}}{\partial \boldsymbol{\Phi}_{n+1}}
    -
    \displaystyle \frac{\partial \mathbf{R}^{3D}}{\partial \mathbf{w}_{n+1}} 
    \Big(
    \displaystyle \frac{\partial \mathbf{R}^{0D}}{\partial \mathbf{w}_{n+1}}
    \Big)^{-1} 
    \displaystyle \frac{\partial \mathbf{R}^{0D}}{\partial \boldsymbol{\Phi}_{n+1}}
    \end{bmatrix}^{(k)}_{n+1}
    \begin{bmatrix}
    \Delta \mathbf{w}^{(k)}_{n+1}
    \\
    \Delta \boldsymbol{\Phi}^{(k)}_{n+1}
    \end{bmatrix}
    \\
    =
    -
    \begin{bmatrix}
    \mathbf{R}^{0D}
    \\
    \mathbf{R}^{3D}
    -
    \displaystyle \frac{\partial \mathbf{R}^{3D}}{\partial \mathbf{w}_{n+1}} 
    \Big(
    \displaystyle \frac{\partial \mathbf{R}^{0D}}{\partial \mathbf{w}_{n+1}}
    \Big)^{-1} \mathbf{R}^{0D}
    \end{bmatrix}^{(k)}_{n+1}
    .
\end{multline}
From this, $\Delta \boldsymbol{\Phi}_{n+1}^{(k)}$ can be determined by solving the linear system from the bottom row
\begin{multline}
    \begin{bmatrix}
    \displaystyle \frac{\partial \mathbf{R}^{3D}}{\partial \boldsymbol{\Phi}_{n+1}}
    -
    \displaystyle \frac{\partial \mathbf{R}^{3D}}{\partial \mathbf{w}_{n+1}} 
    \Big(
    \displaystyle \frac{\partial \mathbf{R}^{0D}}{\partial \mathbf{w}_{n+1}}
    \Big)^{-1} 
    \displaystyle \frac{\partial \mathbf{R}^{0D}}{\partial \boldsymbol{\Phi}_{n+1}}
    \end{bmatrix}^{(k)}_{n+1}
    \Delta \boldsymbol{\Phi}_{n+1}^{(k)}
    =
    -
    \begin{bmatrix}
    \mathbf{R}^{3D}
    -
    \displaystyle \frac{\partial \mathbf{R}^{3D}}{\partial \mathbf{w}_{n+1}} 
    \Big(
    \displaystyle \frac{\partial \mathbf{R}^{0D}}{\partial \mathbf{w}_{n+1}}
    \Big)^{-1} \mathbf{R}^{0D}
    \end{bmatrix}^{(k)}_{n+1}.
\end{multline}
This is identical to the linear system for the uncoupled 3D problem Eq. \eqref{eq:3D_linear_system}, except for additional contributions to the 3D model's residual and tangent from the 0D model. For convenience, we denote the 0D contribution to the 3D tangent $\mathbf{K}^{3D/0D}$, where
\begin{equation} \label{eq:K3D0D}
    \begin{bmatrix}
        \mathbf{K}^{3D/0D}
    \end{bmatrix}
    _{n+1}^{(k)}
    = 
    -
    \begin{bmatrix}
    \displaystyle \frac{\partial \mathbf{R}^{3D}}{\partial \mathbf{w}_{n+1}} 
    \Big(
    \displaystyle \frac{\partial \mathbf{R}^{0D}}{\partial \mathbf{w}_{n+1}}
    \Big)^{-1} 
    \displaystyle \frac{\partial \mathbf{R}^{0D}}{\partial \boldsymbol{\Phi}_{n+1}}
     \end{bmatrix}
      _{n+1}^{(k)}
      .
\end{equation}
As will be shown, rather than considering the 0D contribution to the 3D residual, it is more convenient to consider the entire 0D-modified 3D residual $\mathbf{R}^{3D/0D}$, where
\begin{equation} \label{eq:R3D0D}
    \begin{bmatrix}
    \mathbf{R}^{3D/0D} 
    \end{bmatrix}
    _{n+1}^{(k)}
    = 
    \begin{bmatrix}
    \mathbf{R}^{3D}
    -
    \displaystyle \frac{\partial \mathbf{R}^{3D}}{\partial \mathbf{w}_{n+1}} 
    \Big(
    \displaystyle \frac{\partial \mathbf{R}^{0D}}{\partial \mathbf{w}_{n+1}}
    \Big)^{-1} \mathbf{R}^{0D}
    \end{bmatrix}
    _{n+1}^{(k)}
    .
\end{equation}
Thus, the solution strategy is as follows:
\begin{enumerate}
    \item Approximate 
    $\begin{bmatrix}
    \mathbf{R}^{3D/0D} 
    \end{bmatrix}
    _{n+1}^{(k)}$
    .
    \item Approximate 
    $\begin{bmatrix}
        \mathbf{K}^{3D/0D}
    \end{bmatrix}
    _{n+1}^{(k)}$
    .
    \item Solve the modified 3D linear system
    \begin{equation} \label{eq:modified_3D_linear_system}
        \begin{bmatrix}
    \displaystyle \frac{\partial \mathbf{R}^{3D}}{\partial \boldsymbol{\Phi}_{n+1}}
    +
    \mathbf{K}^{3D/0D}
    \end{bmatrix}^{(k)}_{n+1}
    \\
    \Delta \boldsymbol{\Phi}_{n+1}^{(k)}
    =
    -
    \begin{bmatrix}
    \mathbf{R}^{3D/0D}
    \end{bmatrix}^{(k)}_{n+1}
    .
    \end{equation}
\end{enumerate}
 We then perform the Newton update with $\Delta \boldsymbol{\Phi}_{n+1}^{(k)}$, proceed to the next Newton iteration $k +1$, and repeat until convergence. Note that because we use RK4 (an explicit scheme) for the 0D system, the 0D system does not depend on an updated guess for $\mathbf{w}_{n+1}$, and thus we do not need to compute $\Delta \mathbf{w}_{n+1}^{(k)}$. 

The approximations for 
$\begin{bmatrix}
    \mathbf{R}^{3D/0D} 
    \end{bmatrix}
    _{n+1}^{(k)}$ 
    and 
$\begin{bmatrix}
        \mathbf{K}^{3D/0D}
    \end{bmatrix}
    _{n+1}^{(k)}$ 
are performed using a fixed point iteration operator $F^{0D}(\mathbf{w}_{n+1},\boldsymbol{\Phi}_{n+1})$ for the 0D system, which is introduced next. Then, explicit expressions for the two terms are provided.

\subsubsection{0D fixed point iteration operator}
Here, we introduce the 0D fixed point iteration operator, which is necessary for the 0D residual and tangent approximations to follow. Assume we have an operator $F^{0D}(\mathbf{w}_{n+1},\boldsymbol{\Phi}_{n+1})$ such that the iteration $\mathbf{w}_{n+1}^{(m+1)} = F^{0D}(\mathbf{w}_{n+1}^{(m)},\boldsymbol{\Phi}_{n+1})$ converges to the solution of $\mathbf{R}^{0D}(\mathbf{w}_{n+1},\boldsymbol{\Phi}_{n+1}) = \mathbf{0}$ (for fixed $\boldsymbol{\Phi}_{n+1}$). One can identify this operator for nearly all conceivable 0D solvers, implicit or explicit. In this work, RK4 is used to integrate the 0D system, and the fixed point operator corresponding to RK4, $F^{0D, RK4}(\mathbf{w}_{n+1},\boldsymbol{\Phi}_{n+1})$, is provided in Appendix \ref{sect:RK4_scheme}. 

As was done for $\mathbf{R}^{0D}$ Eq. \eqref{eq:coupled_residuals}, it is convenient to view $F^{0D}$ as a function of $\boldsymbol{\Phi}_{n+1}$ through the coupling flow rates $\mathbf{Q}_{n+1}$ as,
\begin{equation}
    F^{0D}(\mathbf{w}_{n+1}, \boldsymbol{\Phi}_{n+1}) = F^{0D}(\mathbf{w}_{n+1}, \mathbf{Q}_{n+1}(\boldsymbol{\Phi}_{n+1})).
\end{equation}

We briefly list some relevant features of $F^{0D, RK4}$. Because RK4 is an explicit scheme, $F^{0D, RK4}$ is a function of $\boldsymbol{\Phi}_{n+1}$, but is not a function of $\mathbf{w}_{n+1}$. However, for generality and in case one is interested in using an implicit scheme, we retain its dependence on $\mathbf{w}_{n+1}$ in the remainder of the derivation. It is a fixed point operator corresponding to Newton's method, which converges in one iteration; in other words, $\mathbf{w}_{n+1}^{(m+1)} = F^{0D, RK4}(\mathbf{w}_{n+1}^{(m)}, \mathbf{Q}_{n+1}(\boldsymbol{\Phi}_{n+1}))$ converges to the solution of $\mathbf{R}^{0D}(\mathbf{w}_{n+1},\mathbf{Q}_{n+1}(\boldsymbol{\Phi}_{n+1})) = \mathbf{0}$ in only one step regardless of $\mathbf{w}_{n+1}^{(m)}$. Finally, on a practical note, we emphasize that $F^{0D}$ is also a function of $\mathbf{w}_n$ and $\mathbf{Q}_n$, which are assumed to be known. See Appendix \ref{sect:RK4_scheme} for more details.

With $F^{0D}$  described, we continue by deriving explicit expressions for 
 $\begin{bmatrix}
\mathbf{R}^{3D/0D} 
\end{bmatrix}
_{n+1}^{(k)}$ 
and 
$\begin{bmatrix}
\mathbf{K}^{3D/0D}
\end{bmatrix}
_{n+1}^{(k)}$ . 

\subsubsection{0D-modified 3D residual}
In this section, an explicit expression is derived for
$\begin{bmatrix}
    \mathbf{R}^{3D/0D} 
    \end{bmatrix}
    _{n+1}^{(k)}$. 
Recall Eq. \eqref{eq:R3D0D},
\begin{equation*}
 \begin{bmatrix}
    \mathbf{R}^{3D/0D} 
    \end{bmatrix}
    _{n+1}^{(k)}
    = 
    \begin{bmatrix}
    \mathbf{R}^{3D}
    -
    \displaystyle \frac{\partial \mathbf{R}^{3D}}{\partial \mathbf{w}_{n+1}} 
    \Big(
    \displaystyle \frac{\partial \mathbf{R}^{0D}}{\partial \mathbf{w}_{n+1}}
    \Big)^{-1} \mathbf{R}^{0D}
    \end{bmatrix}
    _{n+1}^{(k)}
    .
\end{equation*}
First, use a finite difference approximation as suggested in \cite{matthies2003partitioned},
\begin{equation} \label{eq:R3D_with_w_tilde}
    \begin{bmatrix}
    \mathbf{R}^{3D/0D} 
    \end{bmatrix}
    _{n+1}^{(k)}
     = 
    \begin{bmatrix}
    \mathbf{R}^{3D}
    -
    \displaystyle \frac{\partial \mathbf{R}^{3D}}{\partial \mathbf{w}_{n+1}} 
    \Big(
    \displaystyle \frac{\partial \mathbf{R}^{0D}}{\partial \mathbf{w}_{n+1}}
    \Big)^{-1} \mathbf{R}^{0D}
    \end{bmatrix}
    _{n+1}^{(k)}
    \approx 
    \mathbf{R}^{3D} (\boldsymbol{\Phi}_{n+1}^{(k)}, \tilde{\mathbf{w}}_{n+1}^{(k)})
    ,
\end{equation}
where 
\begin{equation} \label{eq:w_tilde}
   \tilde{\mathbf{w}}_{n+1}^{(k)} 
   = 
   \mathbf{w}_{n+1}^{(k)} 
   - \begin{bmatrix}
       \Big(
       \displaystyle \frac{\partial \mathbf{R}^{0D}}{\partial \mathbf{w}_{n+1}} 
       \Big)^{-1} 
       \mathbf{R}^{0D}
   \end{bmatrix}_{n+1}^{(k)}
    .
\end{equation}
Note that in our case, the 3D residual $\mathbf{R}^{3D}$ is linear in the 0D variables $\mathbf{w}_{n+1}$, so this approximation is exact. Eq. \eqref{eq:w_tilde} is in fact one Newton iteration to solve the 0D system at fixed $\boldsymbol{\Phi}_{n+1}^{(k)}$ (or $\mathbf{Q}_{n+1}^{(k)}$). Thus, we may approximate it using the 0D fixed point operator
\begin{equation} \label{eq:w_tilde_1}
    \tilde{\mathbf{w}}_{n+1}^{(k)} \approx F^{0D}(\mathbf{w}_{n+1}^{(k)}, \boldsymbol{\Phi}_{n+1}^{(k)})
    .
\end{equation}
If $F^{0D}$ is a Newton iteration, as in our case, this approximation is exact. 

In terms of implementation, at each Newton iteration, first compute $\tilde{\mathbf{w}}_{n+1}^{(k)}$ using Eq. \eqref{eq:w_tilde_1}. Recalling $\mathbf{R}^{3D}$ is a function of $\mathbf{w}_{n+1}$ through coupling pressures (Eq. \eqref{eq:coupled_residuals}), next compute modified pressures defined as 
\begin{equation}
    \tilde{\mathbf{P}}_{n+1}^{(k)} = \mathbb{P} \tilde{\mathbf{w}}_{n+1}^{(k)}.
\end{equation}
Finally, evaluate 
\begin{equation} \label{eq:0D_residual_contribution_abstract}
    \begin{bmatrix}
    \mathbf{R}^{3D/0D} 
    \end{bmatrix}
    _{n+1}^{(k)}
    =
    \mathbf{R}^{3D} (\boldsymbol{\Phi}_{n+1}^{(k)}, \tilde{\mathbf{P}}_{n+1}^{(k)})
    .
\end{equation} 
This last step can be implemented as follows in a 3D finite element solver. Splitting the residual into a term from the coupled Neumann boundaries and terms from all other residual contributions (internal stresses, other boundary conditions, etc.), $\begin{bmatrix}
    \mathbf{R}^{3D/0D} 
    \end{bmatrix}
    _{n+1}^{(k)}$ is computed as follows:
\begin{equation} \label{eq:0D_residual_contribution}
    \Big(
    R^{3D/0D}
    \Big)_{n+1, Ai}^{(k)}
    = 
    \text{Uncoupled residual terms}
    + \sum_{m=1}^{n^{cBC}}\int_{\Gamma_{h_c}^{(m)}} N_A \tilde{P}_{n+1, m}^{(k)} n_i d\Gamma
    ,
\end{equation}
where $A$ is the node index, $i$ indexes the spatial dimension, $m$ indexes the coupled Neumann boundaries, of which there are $n^{cBC}$, and $\Gamma_{h_c}^{(m)}$ is the surface corresponding to coupled Neumann boundary $m$. $N_A$ is the shape function for node A, $\tilde{P}_{n+1, m}^{(k)}$ is the $m$th component of $ \tilde{\mathbf{P}}_{n+1}^{(k)}$, and $n_i$ is the $i$th component of the outward surface normal. The integral expression in Eq. \eqref{eq:0D_residual_contribution} is the same as for any other (uncoupled) pressure boundary condition; the only difference is that the value of the pressure $\tilde{P}_{n+1, j}^{(k)}$ is obtained by communicating with the 0D solver. If $\mathbf{R}^{3D}$ contains momentum and continuity components (as in the fluid system Eq. \eqref{eq:fluid_linear_system}), the contribution of the coupled Neumann boundary conditions should be assembled into the \textit{momentum} equation residual.
\\

\noindent \textbf{Remark: } We point out the minor difference in Eq. \eqref{eq:0D_residual_contribution} when considering a 3D fluid - 0D fluid problem vs. a 3D structure - 0D fluid problem. For a 3D fluid, we may consider the fluid domain to be fixed (non-deforming), so the integral is taken over the coupled surface in the reference configuration. For a 3D structure, we typically assume a ``follower pressure load". Thus, if the structure deforms, the integral is taken over the coupled surface in the current (deformed) configuration with the current surface normal vector, corresponding to timestep $n+1$.

\subsubsection{0D contribution to 3D tangent} \label{sect:0D_tangent_contribution}
In this section, an explicit expression is derived for 
$\begin{bmatrix}
        \mathbf{K}^{3D/0D}
\end{bmatrix}_{n+1}^{(k)}$
. 
Recall Eq. \eqref{eq:K3D0D},
\begin{equation*}
    \begin{bmatrix}
        \mathbf{K}^{3D/0D}
    \end{bmatrix}
    _{n+1}^{(k)}
    = 
    -
    \begin{bmatrix}
    \displaystyle \frac{\partial \mathbf{R}^{3D}}{\partial \mathbf{w}_{n+1}} 
    \Big(
    \displaystyle \frac{\partial \mathbf{R}^{0D}}{\partial \mathbf{w}_{n+1}}
    \Big)^{-1} 
    \displaystyle \frac{\partial \mathbf{R}^{0D}}{\partial \boldsymbol{\Phi}_{n+1}}
     \end{bmatrix}
      _{n+1}^{(k)}
      .
\end{equation*}
Following \cite{chan1985approximate}, define a matrix $\mathbf{C}$
\begin{equation}
    \mathbf{C} 
    = 
    \Big(
    \frac{\partial \mathbf{R}^{0D}}{\partial \mathbf{w}_{n+1}}
    \Big)^{-1} 
    \frac{\partial \mathbf{R}^{0D}}{\partial \boldsymbol{\Phi}_{n+1}},
\end{equation}
so that
\begin{equation*}
    \begin{bmatrix}
        \mathbf{K}^{3D/0D}
    \end{bmatrix}
    _{n+1}^{(k)}
    =
     -
    \begin{bmatrix}
    \displaystyle \frac{\partial \mathbf{R}^{3D}}{\partial \mathbf{w}_{n+1}} 
    \mathbf{C}
     \end{bmatrix}
      _{n+1}^{(k)}.
\end{equation*}
The key approximation is provided in \cite{chan1985approximate}, in which it was shown that $\mathbf{C}$ can be reasonably approximated by
\begin{equation}
    \mathbf{C} 
    \approx 
    - \frac{\partial F^{0D}(\mathbf{w}_{n+1}, \boldsymbol{\Phi}_{n+1})}{\partial \boldsymbol{\Phi}_{n+1}}
    =
    - \frac{\partial F^{0D}(\mathbf{w}_{n+1}, \mathbf{Q}_{n+1})}{\partial \mathbf{Q}_{n+1}}
    \frac{\partial \mathbf{Q}_{n+1}}{\partial \mathbf{U}_{n+1}}
    \frac{\partial \mathbf{U}_{n+1}}{\partial \boldsymbol{\Phi}_{n+1}}
    ,
\end{equation}
where we have used the chain rule and the fact that $F^{0D}$ is a function of $\boldsymbol{\Phi}_{n+1}$ through the flow rates $\mathbf{Q}_{n+1}$ and the nodal velocities $\mathbf{U}_{n+1}$ (Eqs. \eqref{eq:Q3D} and \eqref{eq:genalpha_U}).

For the other term, $\displaystyle \frac{\partial \mathbf{R}^{3D}}{\partial \mathbf{w}_{n+1}}$, we also apply the chain rule and the fact that $\mathbf{R}^{3D}$ is a function of $\mathbf{w}_{n+1}$ through the pressures $\mathbf{P}_{n+1}$ (Eq. \eqref{eq:P_relation})

\begin{equation}
    \frac{\partial \mathbf{R}^{3D}}{\partial \mathbf{w}_{n+1}} = \frac{\partial \mathbf{R}^{3D}}{\partial \mathbf{P}_{n+1} } \frac{d \mathbf{P}_{n+1} }{d \mathbf{w}_{n+1}} = \frac{\partial \mathbf{R}^{3D}}{\partial \mathbf{P}_{n+1} } \mathbb{P}
    .
\end{equation}
Thus, the 0D tangent matrix contribution is 
\begin{equation} \label{eq:0D_tangent_contribution_abstract}
    \begin{bmatrix}
        \mathbf{K}^{3D/0D}
    \end{bmatrix}
    _{n+1}^{(k)}
    =
    \begin{bmatrix}
    \displaystyle \frac{\partial \mathbf{R}^{3D}}{\partial \mathbf{P}_{n+1} } 
    \displaystyle \mathbb{P}
    \displaystyle \frac{\partial F^{0D}(\mathbf{w}_{n+1}, \mathbf{Q}_{n+1})}{\partial \mathbf{Q}_{n+1}}
    \displaystyle \frac{\partial \mathbf{Q}_{n+1}}{\partial \mathbf{U}_{n+1}}
    \displaystyle \frac{\partial \mathbf{U}_{n+1}}{\partial \boldsymbol{\Phi}_{n+1}}
    \end{bmatrix}
    _{n+1}^{(k)}.
\end{equation}
All partial derivatives are computed analytically (see Appendix \ref{sect:0D_tangent_contribution_derivation} for details) except for 
\begin{equation*}
    \begin{bmatrix}
\mathbb{P} \displaystyle \frac{\partial F^{0D}(\mathbf{w}_{n+1}, \mathbf{Q}_{n+1})}{\partial \mathbf{Q}_{n+1}}
\end{bmatrix}_{n+1}^{(k)}.
\end{equation*}
This term, which is denoted by $\mathbf{M}$, is computed in a finite difference manner by communicating with the 0D solver as
\begin{align*}
    M_{ij}
    &=
    \begin{bmatrix}
    \mathbb{P}_{ip}
    \Big( 
    \displaystyle \frac{\partial F^{0D}(\mathbf{w}_{n+1}, \mathbf{Q}_{n+1})}{\partial \mathbf{Q}_{n+1}}
    \Big)_{pj}
    \end{bmatrix}_{n+1}^{(k)}
    \\
    &\approx
    \frac{\mathbb{P}_{ip} F^{0D}_p(\mathbf{w}_{n+1}^{(k)}, \mathbf{Q}_{n+1}^{(k)} + \epsilon \mathbf{e}_j) - \mathbb{P}_{ip} F^{0D}_p(\mathbf{w}_{n+1}^{(k)}, \mathbf{Q}_{n+1}^{(k)})}{\epsilon},
\end{align*}
where $\mathbf{e}_j$ is the $j$th unit vector and $\epsilon$ is a small numerical perturbation.
Note that for the evaluation of
$\begin{bmatrix}
\mathbf{R}^{3D/0D} 
\end{bmatrix}
_{n+1}^{(k)}$,
we already require
\begin{equation*}
    \mathbb{P} F^{0D}(\mathbf{w}_{n+1}^{(k)}, \mathbf{Q}_{n+1}^{(k)}) 
    = \mathbb{P} \tilde{\mathbf{w}}_{n+1}^{(k)} =  \tilde{\mathbf{P}}_{n+1}^{(k)}.
\end{equation*}
The slightly perturbed quantity, $\mathbb{P} F^{0D}(\mathbf{w}_{n+1}^{(k)}, \mathbf{Q}_{n+1}^{(k)} + \epsilon \mathbf{e}_j)$, can be computed in exactly the same manner. We define $\tilde{\mathbf{w}}_{n+1, \epsilon_j}^{(k)}$ and $\tilde{\mathbf{P}}_{n+1, \epsilon_j}^{(k)}$ such that
\begin{equation*}
    \mathbb{P} F^{0D}(\mathbf{w}_{n+1}^{(k)}, \mathbf{Q}_{n+1}^{(k)} + \epsilon \mathbf{e}_j) 
    = \mathbb{P} \tilde{\mathbf{w}}_{n+1, \epsilon_j}^{(k)} =  \tilde{\mathbf{P}}_{n+1, \epsilon_j}^{(k)}.
\end{equation*}
Thus, we may write
\begin{equation} \label{eq:M_matrix}
    M_{ij}
    = 
    \frac{
    \tilde{P}_{n+1,\epsilon_j, i}^{(k)}
    - 
    \tilde{P}_{n+1,i}^{(k)}
    }
    {\epsilon}
    .
\end{equation}
In this form, it is revealed that $M_{ij}$ is a resistance-like quantity that describes how pressure at coupled surface $i$ changes with flow rate at coupled surface $j$. As argued in \cite{moghadam2013modular}, the off-diagonal entries of $M_{ij}$ are ignored in this work. 

The resulting 0D tangent contribution is given as (again, see Appendix \ref{sect:0D_tangent_contribution_derivation} for derivation) 
\begin{equation} \label{eq:0D_tangent_contribution}
    \Big(
    K^{3D/0D}
    \Big)_{n+1, AiBj}^{(k)} 
    =
    \sum_{l=1}^{n^{cBC}} \sum_{m=1}^{n^{cBC}} \gamma \Delta t M_{lm}
    \int_{\Gamma_{h_c}^{(l)}} N_A n_i d\Gamma
    \int_{\Gamma_{h_c}^{(m)}} N_B n_j d\Gamma
    .
\end{equation}
The variable definitions are the same as for $\mathbf{R}^{3D/0D}$ (a vector), with the addition of indices ($B$, $j$, $m$) since we are now dealing with a matrix. 
Analogous to the 0D-modified 3D residual, the 0D tangent contribution should be assembled into the \textit{momentum-acceleration} block of the tangent matrix, if applicable.
\\

\noindent \textbf{Remark: } For a 3D structure - 0D fluid problem, the integrals in Eq. \eqref{eq:0D_tangent_contribution} should be taken over the coupled surface in the current (deformed) configuration with the current surface normal vector, corresponding to timestep $n+1$.
\\

\noindent \textbf{Remark: } Eq. \eqref{eq:0D_residual_contribution} (0D-modified 3D residual) and Eq. \eqref{eq:0D_tangent_contribution} (0D contribution to the 3D tangent) can be found in \cite{moghadam2013modular, esmaily2013new}, in slightly different notation. However, we obtain these expressions through a new derivation involving ANM and show that they apply not only to the 3D fluid - 0D fluid problem, but also to the 3D structure - 0D fluid problem.

\subsubsection{Summary} \label{sect:summary}
\begin{figure}
    \centering
    \includegraphics[width = 1.0\textwidth]{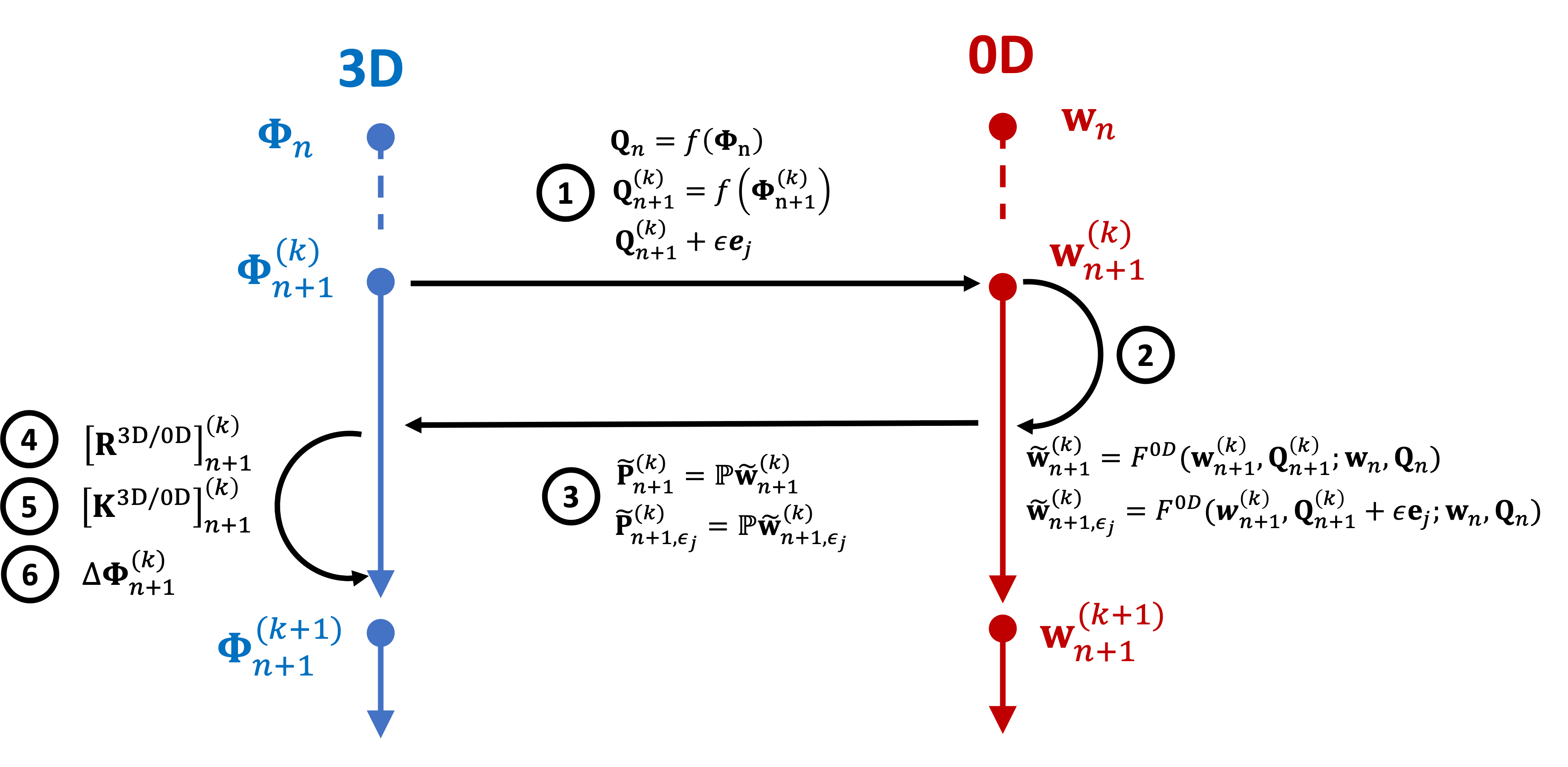}
    \caption{Communication diagram between 3D and 0D solvers. At each timestep $n$, we iterate until convergence. At each Newton iteration $k$, we perform 6 computations/communications. These steps correspond to the coupling algorithm summary in Section \ref{sect:summary}.}
    \label{fig:communication_pattern}
\end{figure}

The coupling strategy is summarized below. See Fig. \ref{fig:communication_pattern} for a graphic representation of the required computations and communications. At each Newton iteration $k$, 
\begin{enumerate}
    \item Compute $\mathbf{Q}^{(k)}_{n+1} = f(\boldsymbol{\Phi}_{n+1}^{(k)})$, where $f$ is some function according to Eqs. \eqref{eq:Q3D} and \eqref{eq:genalpha_U}, and communicate with the 0D solver. The 3D solver must also compute and send $\mathbf{Q}_{n} = f(\boldsymbol{\Phi}_n)$, which is also required by the 0D solver (see Appendix \ref{sect:RK4_scheme} for more details).
    \item Compute 
    \begin{align*}
        \tilde{\mathbf{w}}_{n+1}^{(k)} 
        &= 
        F^{0D}(
        \mathbf{w}_{n+1}^{(k)}, \mathbf{Q}_{n+1}^{(k)}; \mathbf{w}_{n}, \mathbf{Q}_{n}
        ),
        \\
        \tilde{\mathbf{w}}_{n+1, \epsilon_j}^{(k)} 
        &= 
        F^{0D}(
        \mathbf{w}_{n+1}^{(k)}, \mathbf{Q}_{n+1}^{(k)} + \epsilon \mathbf{e}_j; \mathbf{w}_{n}, \mathbf{Q}_{n}
        ), 
        \hspace{10 pt} j \in \{1, \ldots, n_{bc} \}.
    \end{align*}
    This requires $n_{bc} + 1$ calls to the 0D solver.
    
    \item Compute 
    \begin{align*}
        \tilde{\mathbf{P}}_{n+1}^{(k)} 
        &= 
        \mathbb{P} \tilde{\mathbf{w}}_{n+1}^{(k)},
        \\
        \tilde{\mathbf{P}}_{n+1, \epsilon_j}^{(k)} 
        &= 
        \mathbb{P} \tilde{\mathbf{w}}_{n+1, \epsilon_j}^{(k)}, 
        \hspace{10 pt}
        j \in \{1, \ldots, n_{bc} \}.
    \end{align*} 
    This requires the 0D to simply extract the proper components of $\tilde{\mathbf{w}}_{n+1}^{(k)}$ and $\tilde{\mathbf{w}}_{n+1, \epsilon_j}^{(k)}$ and send them to 3D.
    
    \item Given 
    $\tilde{\mathbf{P}}_{n+1}^{(k)} $, 
    compute 
    $\begin{bmatrix}
    \mathbf{R}^{3D/0D} 
    \end{bmatrix}
    _{n+1}^{(k)}$ 
    using Eq. \eqref{eq:0D_residual_contribution}.
    
    \item Given 
    $\tilde{\mathbf{P}}_{n+1}^{(k)} $ 
    and 
    $\tilde{\mathbf{P}}_{n+1, \epsilon_j}^{(k)}$, $j \in \{1, \ldots, n_{bc} \}$, 
    compute $\mathbf{M}$ using Eq. \eqref{eq:M_matrix} and construct 
    $\begin{bmatrix}
        \mathbf{K}^{3D/0D}
    \end{bmatrix}_{n+1}^{(k)}$
    using Eq. \eqref{eq:0D_tangent_contribution}.
    \item Solve the modified 3D linear system Eq. \eqref{eq:modified_3D_linear_system} 
    \begin{equation*}
        \begin{bmatrix}
        \displaystyle \frac{\partial \mathbf{R}^{3D}}{\partial \boldsymbol{\Phi}_{n+1}}
        +
        \mathbf{K}^{3D/0D}
        \end{bmatrix}^{(k)}_{n+1}
        \\
        \Delta \boldsymbol{\Phi}_{n+1}^{(k)}
        =
        -
        \begin{bmatrix}
        \mathbf{R}^{3D/0D}
        \end{bmatrix}^{(k)}_{n+1},
    \end{equation*}
    for $\Delta \boldsymbol{\Phi}_{n+1}^{(k)}$.
\end{enumerate}

The 3D degrees of freedom are updated with $\Delta \boldsymbol{\Phi}_{n+1}^{(k)}$, the Newton iteration index is incremented $k \rightarrow k+1$, and the process repeats until $\begin{bmatrix}
\mathbf{R}^{3D/0D}
\end{bmatrix}^{(k)}_{n+1}$ 
falls below a prescribed relative or absolute tolerance.

\subsection{Numerical considerations}
The 0D tangent contribution 
$
\begin{bmatrix}
\mathbf{K}^{3D/0D}
\end{bmatrix}
_{n+1}^{(k)}
$
is dense and adding it explicitly to 
$
\begin{bmatrix}
\displaystyle \frac{\partial \mathbf{R}^{3D}}{\partial \boldsymbol{\Phi}_{n+1}}
\end{bmatrix}
_{n+1}^{(k)}
$
may deteriorate the performance of the linear solver, which typically takes advantage of the sparse structure of tangent matrices arising in finite element solvers \cite{evans1973analysis}. To address this, we rewrite Eq. \eqref{eq:0D_tangent_contribution} as
\begin{equation}
    \Big(
    K^{3D/0D}
    \Big)_{n+1, AiBj}^{(k)} 
    = 
    \sum_{l=1}^{n^{cBC}} \sum_{m=1}^{n^{cBC}} \gamma \Delta t M_{lm} (S_l)_{Ai}  (S_m)_{Bj}, \hspace{10pt} (S_l)_{Ai} = \int_{\Gamma_{h_c}^{(l)}} N_A n_i d\Gamma
    .
\end{equation}
That is, 
$\begin{bmatrix}
\mathbf{K}^{3D/0D}
\end{bmatrix}
_{n+1}^{(k)}$ 
is the sum of rank 1 matrices. It is more efficient to store the vector $\mathbf{S}_k$ separately and apply when needed (in matrix multiplication and in preconditioning) than to explicitly form the outer product and add it to the tangent \cite{esmaily2013new}.

The 0D tangent contribution also increases the condition number of the linear system, proportional to the resistance of the coupled Neumann boundaries (i.e., $\mathbf{M}$). This can cause poor performance of standard iterative linear solvers. Resistance-based preconditioning is an effective remedy \cite{seo2019performance, esmaily2013new}.

\subsection{Capping non-closed surfaces and consequences} \label{sect:capping}
\begin{figure}
    \centering
    \includegraphics[width=0.6\linewidth]{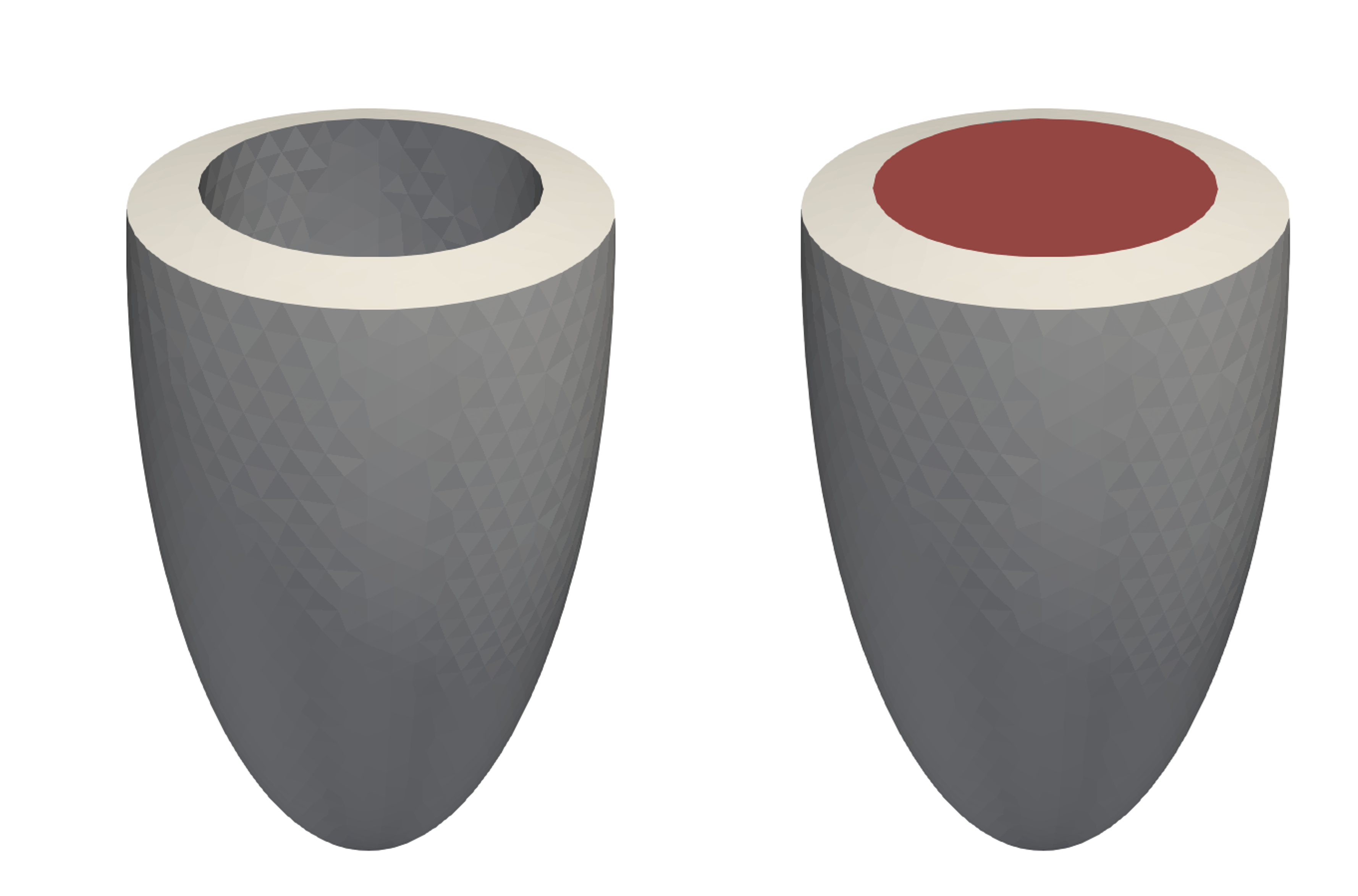}
    \caption{Left: An idealized LV. Right: The same model with a cap surface shown in red, which, combined with the endocardial surface, defines a closed surface with which to compute volume/flow rate.
}
    \label{fig:capped_LV}
\end{figure}

The proposed coupling method requires computing flow rates from the 3D domain. For a 3D structure, this flow rate is identical to the negative rate of change of the enclosed volume. An important complication arises if the structure does not enclose a volume, such as in modeling the mechanics of the heart muscle. When the cardiac valves are closed, each cardiac chamber encloses a volume. However, the valves are often ignored when modeling the heart. As an example, Fig. \ref{fig:capped_LV} (left) shows an idealized LV, where the inflow and outflow valves are omitted. Here, we would like to couple the endocardial surface to a 0D fluid model, but the endocardial surface is not closed. Not accounting for this will lead to an inaccurately computed flow rate. In this work, we address this issue by introducing a ``cap," a surface that closes the endocardial surface, thereby defining an enclosed fluid-tight volume (Fig. \ref{fig:capped_LV} right). Stated formally, if a coupled surface $\Gamma_{h_c}^{(i)}$ is not closed, we consider a cap surface $\Gamma_{h_c, cap}^{(i)}$ such that $\Gamma_{h_c}^{(i)} \cup \Gamma_{h_c, cap}^{(i)}$ is a closed surface. 

Because there is some flexibility in defining the cap surface, in our formulation, we leave the responsibility of constructing it to the user.
Methods exist for constructing such surfaces, such as the ear clipping algorithm \cite{eberly2008triangulation} or the \texttt{vtkFillHolesFilter} from VTK \cite{schroeder2006visualization}, the latter being used in this work.

We note that this cap surface is used only for computing flow rate and should not be treated in the same manner as a boundary of the 3D model. In particular, even though the cap surface is used to compute flow rate, the coupling pressure is not applied to the cap surface. The consequences are as follows. When computing flow rates, the cap surface is included in the integral in Eq. \eqref{eq:Q3D}, 
 \begin{equation*}    
     Q_i 
     = 
     \sum_A 
     \Big(
     \int _{\Gamma_{h_c}^{(i)}} N_A (\mathbf{U})_A \cdot \mathbf{n} d\Gamma 
     + 
     \int _{\Gamma_{h_c, cap}^{(i)}} N_A (\mathbf{U})_A \cdot \mathbf{n} d\Gamma
     \Big)
     .
\end{equation*}
When computing the 0D-modified residual, the cap surface is not included in the pressure integral in Eq. \eqref{eq:0D_residual_contribution}; it is unchanged,
\begin{equation*}
    \Big(
    R^{3D/0D}
    \Big)_{n+1, Ai}^{(k)}
    = 
    \text{Uncoupled residual terms}
    + \sum_{j=1}^{n^{cBC}}\int_{\Gamma_{h_c}^{(j)}} N_A \tilde{P}_{n+1, j}^{(k)} n_i d\Gamma.
\end{equation*}
When computing the 0D tangent contribution, the cap surface is ignored in the first integral, but included in the second integral in Eq. \eqref{eq:0D_tangent_contribution},
\begin{equation*}
    \Big(
    K^{3D/0D}
    \Big)_{n+1, AiBj}^{(k)} 
    =
    \sum_{l=1}^{n^{cBC}} \sum_{m=1}^{n^{cBC}} \gamma \Delta t M_{lm}
    \int_{\Gamma_{h_c}^{(l)}} N_A n_i d\Gamma
    \Big(
    \int_{\Gamma_{h_c}^{(m)}} N_B n_j d\Gamma
    +
    \int_{\Gamma_{h_c, cap}^{(m)}} N_B n_j d\Gamma
    \Big)
    .
\end{equation*}

Note that the expressions above assume the cap surface does not introduce any new nodes in the interior of the surface. If the cap included interior nodes, we would have to interpolate velocities from the boundary nodes (where the velocity is well-defined) to these interior nodes, and the above expressions would have to account for this interpolation. We leave this problem for future work.

\section{Results} \label{sect:results}
The present coupling method is implemented in svFSI, a multiphysics finite element solver for cardiovascular modeling \cite{zhu2022svfsi}. In the following sections, we demonstrate the coupling method in several illustrative and clinically relevant examples, as well as assess the convergence behavior and cost of our method.

\subsection{Idealized left ventricle with open-loop circulation} \label{sect:LV_sim}

\begin{figure}
    \centering
    \includegraphics[width=1.0\columnwidth]{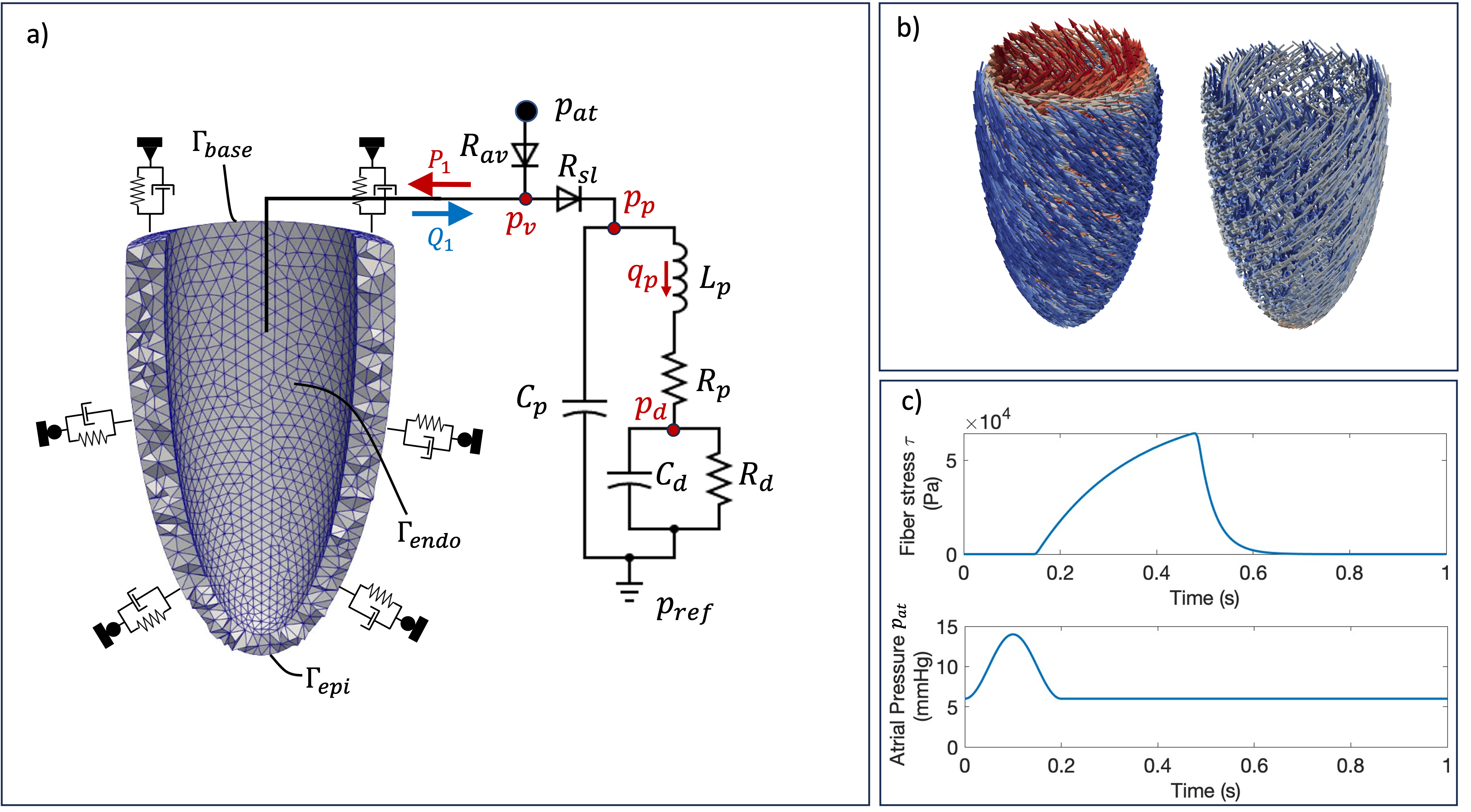}
    \caption{Problem setup for the coupled idealized LV example. a) An idealized LV model, shown in cut-view, is coupled to an open-loop LPN of the systemic circulation (LPN parameters are given in Table \ref{tab:LV_sim_params}). The LV is supported by Robin boundary conditions (denoted by the spring-dashpot assemblies) on the base and epicardial surfaces, and periodic active stress along fiber directions causes the LV to contract.
    b) Fiber $\mathbf{f}$ (left) and sheet $\mathbf{s}$ (right) orientation fields for the idealized LV model. Arrows denote the local fiber or sheet direction, and they are colored by their component along the longitudinal axis. Fiber angles (relative to circumferential) vary from $+60^\circ$ on the endocardial surface to $-60^\circ$ degrees on the epicardial surface.
    c) Fiber stress curve (top) and prescribed atrial pressure curve (bottom). Here, we assume the cardiac cycle duration is 1s. Atrial systole is set to begin at $t_{sys, a} = 0 \unit{ms}$ and lasts for a duration $T_{sys, a} = 200 \unit{ms}$. Atrial pressure ranges between 6 \unit{mmHg} and 14 \unit{mmHg}. Ventricular systole is set to begin at $t_{sys,v} = 143 \unit{ms}$, and the active stress reaches a maximum value of approximately 65 \unit{kPa} (488 \unit{mmHg}). 
}
    \label{fig:LV_sim_setup}
\end{figure}

In this example, we simulate a LV with active contraction, pumping blood through an open-loop LPN model of the systemic circulation. The idealized LV geometry is that of a truncated, prolate ellipsoid \cite{land2015verification} (Fig. \ref{fig:LV_sim_setup}a), with three surfaces: epicardial (outer) $\Gamma_{epi}$, endocardial (inner) $\Gamma_{endo}$, and basal (top) $\Gamma_{base}$.
We also construct a fiber orientation field \cite{land2015verification}, which is used in a fiber-based constitutive model of the myocardium, described next. The fiber orientations $\mathbf{f}$ vary linearly from $+60^\circ$ (relative to circumferential) on the endocardial surface to $-60^\circ$ on the epicardial surface. 
The material model also requires a ``sheet" orientation field $\mathbf{s}$, which in this work is perpendicular to $\mathbf{f}$ and to the ellipsoidal normal direction $\mathbf{n}$. See Fig. \ref{fig:LV_sim_setup}b for visualization.

The myocardium is modeled with a Holzapfel-Ogden strain energy,
plus a quadratic volumetric penalty term, as well as a viscous pseudo-potential and active stress along fiber directions to recapitulate cardiac contraction \cite{holzapfel2009constitutive, pfaller2019importance}. The second Piola-Kirchhoff stress is thus given by
\begin{equation} \label{eq:2PK_stress}
    \mathbf{S} = \frac{\partial}{\partial \mathbf{E}}(\psi_{HO} + \psi_{vol}) + \frac{\partial}{\partial \dot{\mathbf{E}}}(\psi_{visc}) + \mathbf{S}_{act},
\end{equation}
\begin{equation} \label{eq:HO_model}
    \psi_{HO} = \frac{a}{2b}
    \Big( 
    e^{b(\Bar{I}_1 - 3)} - 1
    \Big)
    +
    \frac{a_{fs}}{2b_{fs}}
    \Big(
    e^{b_{fs}I_{8,fs}^2} - 1
    \Big)
    +
    \sum_{i \in \{f,s\}}
    \chi(I_{4,i})
    \frac{a_i}{2b_i}
    \Big(
    e^{b_i(I_{4,i} - 1)^2} - 1
    \Big),
\end{equation}
\begin{equation}
    \psi_{vol} = \frac{\kappa}{2}(1-J)^2,
\end{equation}
\begin{equation}
    \psi_{visc} = \frac{\eta}{2} \text{tr}(\dot{\mathbf{E}}^2),
\end{equation}
where $a_i$ and $b_i$ are material parameters, $\kappa$ is a volumetric penalty parameter, and the strain invariants are defined as
\begin{equation}
\begin{split}
    \Bar{I}_1 &= J^{-2/3} I_1 \text{ , where } I_1 = \text{tr}(\mathbf{C}),
    \\
    I_{4,f} &= \mathbf{f} \cdot \mathbf{C} \mathbf{f},
    \\
    I_{4,s} &= \mathbf{s} \cdot \mathbf{C} \mathbf{s},
    \\
    I_{8, fs} &= \mathbf{f} \cdot \mathbf{C} \mathbf{s}.
\end{split}
\end{equation}
$\eta$ is the viscosity and $\dot{\mathbf{E}}$ is the rate of Green-Lagrange strain tensor. 
$\chi(x)$ is a smoothed Heaviside function centered at $x = 1$ with smoothing parameter $k$
\begin{equation} \label{eq:heaviside}
    \chi(x) = \frac{1}{1 + e^{-k(x-1)}}
\end{equation}
Finally, the active stress is applied along fiber directions to recapitulate cardiac contraction,
\begin{equation} \label{eq:active_stress}
    \mathbf{S}_{act} = \tau(t) \cdot \mathbf{f} \otimes \mathbf{f}.
\end{equation}
$\tau(t)$ is determined using the model from \cite{pfaller2019importance}, yielding the active stress curve shown in Fig. \ref{fig:LV_sim_setup}c (top). 

On $\Gamma_{epi}$, we apply a Robin boundary condition in the normal direction only, following \cite{pfaller2019importance} to mimic the effect of the pericardium. On $\Gamma_{base}$, we apply Robin boundary conditions in all directions. On $\Gamma_{endo}$, we use a coupled Neumann boundary condition (Fig. \ref{fig:LV_sim_setup}a). The initial LV geometry is set as the stress-free reference configuration, and the simulation is initialized with zero displacements and velocities.

In this example, the 0D fluid is an open-loop Windkessel-type model of the systemic circulation \cite{pfaller2019importance}, shown in Fig. \ref{fig:LV_sim_setup}a.  Performing a nodal analysis on this LPN yields the following set of ODEs.

\begin{equation} \label{eq:open_loop_LPN}
\begin{split}
    \frac{p_v-p_{at}}{R_{av}} + \frac{p_v - p_p}{R_{sl}} - Q_1 &= 0, \\
    q_p - \frac{p_v - p_p}{R_{sl}} + C_p \dot{p}_p &= 0,\\
    q_p + \frac{p_d - p_p}{R_p} + \frac{L_p}{R_p}\dot{q}_p &= 0,\\
    \frac{p_d - p_{ref}}{R_d} - q_p + C_d \dot{p}_d &= 0.
\end{split}
\end{equation}
The atrial pressure $p_{at}$ is a prescribed function of time, given by
\begin{equation}
    p_{at} = 
    \begin{cases}
         \frac{\Delta p_{at}}{2} 
        \Big(
        1 - \cos{\frac{2\pi (t-t_{sys,a})}{T_{sys,a}}} \Big) + p_{at0}, &
        \text{for $t_{sys,a}< t < t_{sys,a}+ T_{sys,a}$},
        \\
        p_{at0}, & \text{otherwise}.
    \end{cases} 
\end{equation}
See Fig. \ref{fig:LV_sim_setup}c (bottom) for plot.

The atrioventricular (av) and semilunar (sl) valves are modeled as diodes with nonlinear resistances $R_{av}$ and $R_{sl}$ that depend on the pressure differential on either side as below:
\begin{equation}
    R_{av} = R_{min} + (R_{max} - R_{min}) S^+ (p_v - p_{at}),
\end{equation}
\begin{equation}
    R_{sl} = R_{min} + (R_{max} - R_{min}) S^+ (p_p - p_v),
\end{equation}
where $S^+$ is a sigmoid function with steepness parameter $k_p$.

To cast into the general form in Section \ref{sect:coupled_problem}, we identify the 0D unknowns as $\mathbf{w} = [ p_v, p_p, p_d, q_p]^T$, while the coupled Dirichlet boundary data are $\mathbf{Q} = [Q_1]^T$. The uncoupled Neumann boundary data are $\mathbf{p} = [p_{at}, p_{ref}]^T$, and there are no uncoupled Dirichlet boundary data $\mathbf{q}_u = []$. The 0D variables are initialized using values from \cite{pfaller2019importance}. The coupled problem for this example can be summarized as
\begin{equation} 
    \begin{cases}
    \mathcal{P}^{0D}(\mathbf{w}, t, [\mathbf{q}_{u}, \mathbf{Q}], \mathbf{p}) = \mathbf{0} \hspace{7pt} \text{given by Eq. \eqref{eq:open_loop_LPN}},
    \\
    \text{Initial conditions from \cite{pfaller2019importance}},
    \\
    \\
    \mathcal{P}^{3D,struct}(\mathbf{u}, \mathbf{x}, t) = \mathbf{0},
    \\
    \text{Initialize with zero displacement and velocity},
    \\
    \text{(Uncoupled) Robin boundary conditions on } \Gamma_{epi}
    \text{ and } \Gamma_{base},
    \\
    \boldsymbol{\sigma} \cdot \mathbf{n} = -P_1 \mathbf{n}, \text{ on } \Gamma_{endo},
    \\
    \\
   P_1(t) = w_1(t) = p_v(t),
    \\
    Q_1(t) = \int_{\Gamma_{endo}(t)} \mathbf{u}(t) \cdot \mathbf{n}(t) d\Gamma.
    
    \end{cases}
\end{equation}
\newline
\newline
A complete table of parameters for this simulation can be found in Appendix \ref{sect:sim_params} Table \ref{tab:LV_sim_params}.
\newline

\begin{figure} [h]
    \centering
    \includegraphics[width=0.75\columnwidth]{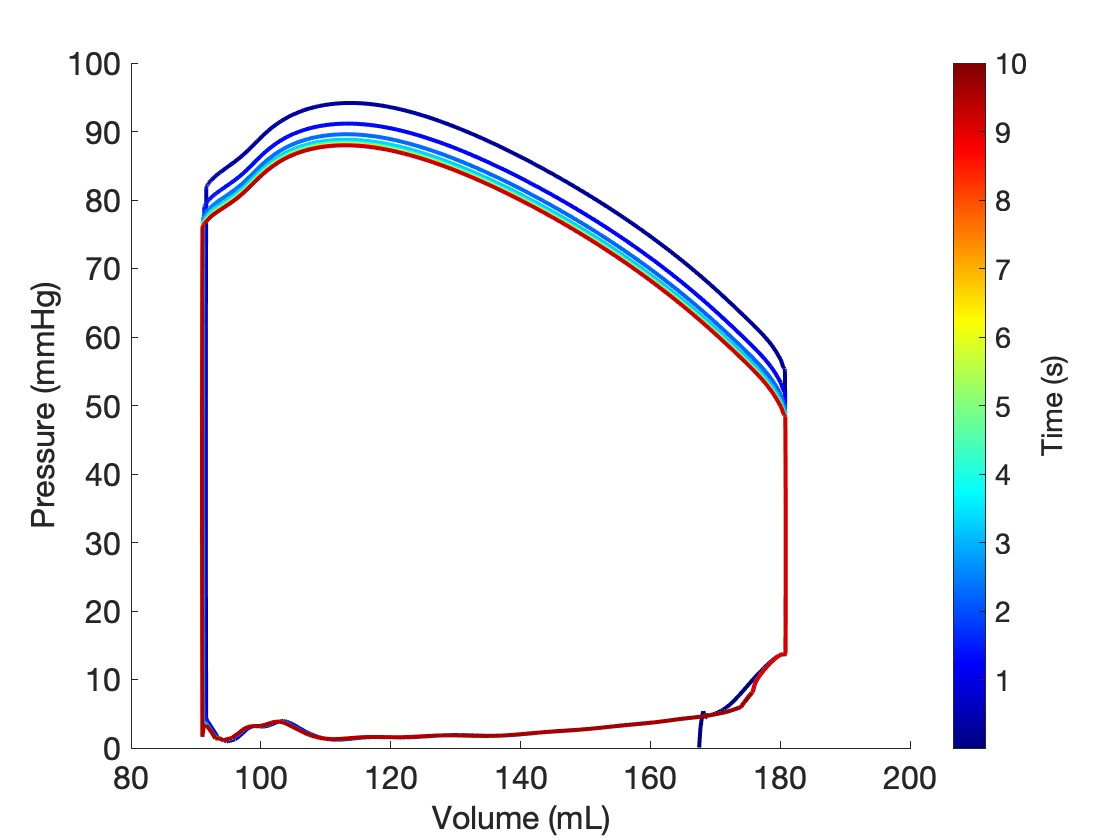}
    \caption{The PV loop for the LV is plotted for 10 cardiac cycles. The PV loop spans a realistic range of pressure and volume,  exhibits a clear distinction of the four cardiac phases, and reaches a limit cycle after about 5 cardiac cycles.
    A movie showing the LV deformation synchronized with the PV loop is provided here: \url{https://drive.google.com/file/d/17ZCVQq8p-EOBrLB3C7olLT2z2uyNTkWu/view?usp=sharing}.
}
    \label{fig:LV_PV_color}
\end{figure}

Fig. \ref{fig:LV_PV_color} shows the pressure-volume (PV) loop over 10 cardiac cycles obtained from this coupled simulation. Although the initial volume of the idealized LV is higher than a physiological heart, the ranges of pressure and volume are in a physiological range. The stroke volume is approximately 90 \unit{mL} and the resulting ejection fraction is 50\%, which is comparable to normal physiological values and values reported in other computational studies with different LV geometries \cite{kosaraju2017left, pfaller2019importance, regazzoni2022cardiac}. There is a clear delineation of the four cardiac phase -- isovolumic contraction, ejection, isovolumic relaxation, and filling. Additionally, the PV loop reaches a limit cycle after about 5 cardiac cycles.

\subsection{Temporal convergence} \label{sect:temporal_convergence}
\begin{figure}
    \centering
    \includegraphics[width = 0.8\textwidth]{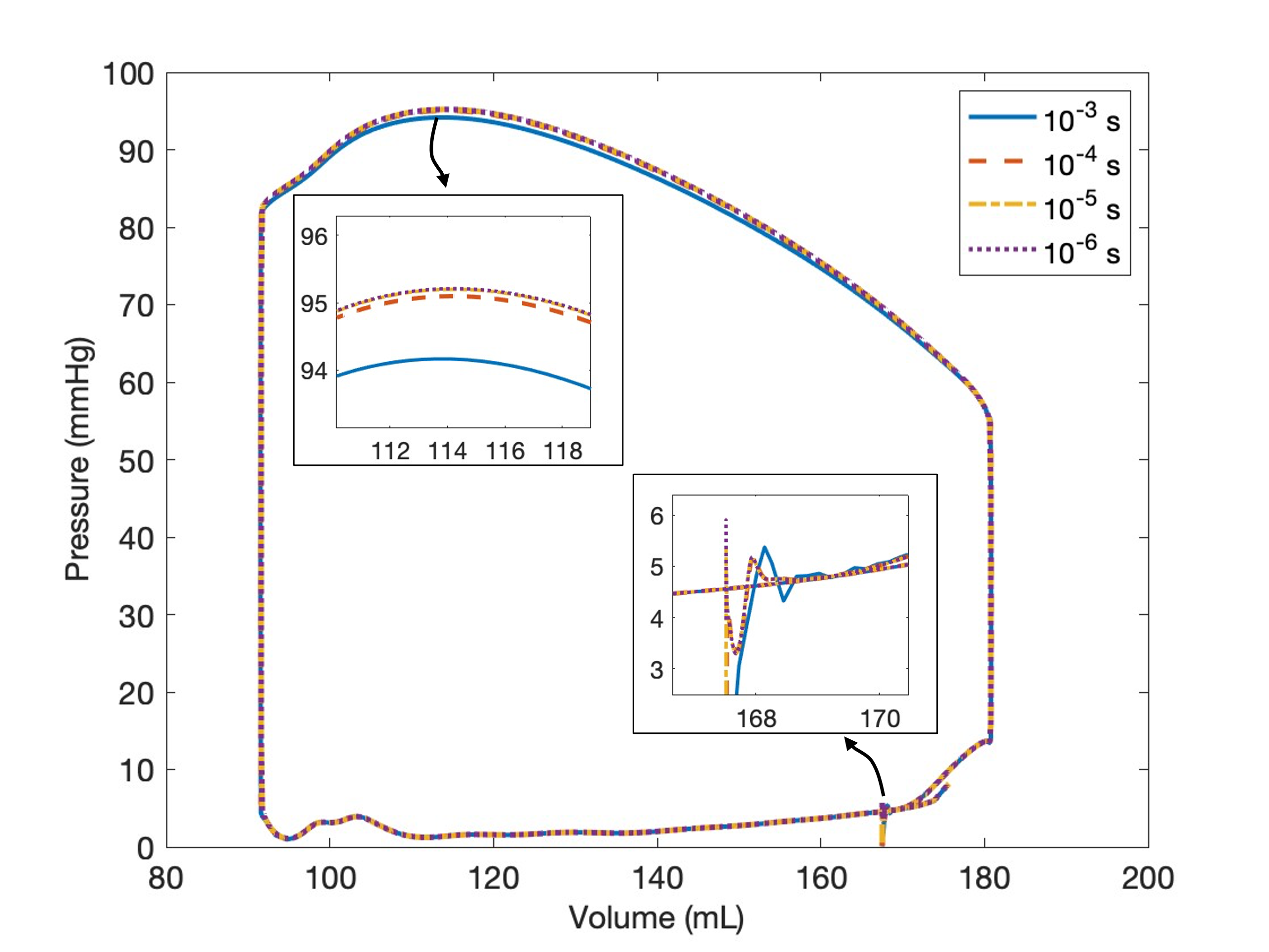}
    \caption{PV loops for a single cardiac cycle of the idealized LV model (Section \ref{sect:LV_sim}) for decreasing timestep size. Insets show zoomed in view at the top and bottom-right portions of the loop. PV loops converge as the timestep size decreases.}
    \label{fig:PV_convergence}
\end{figure}
 
In this section, we present preliminary results on the temporal convergence of the coupling method, applied to the idealized LV model (Section \ref{sect:LV_sim}). In Fig. \ref{fig:PV_convergence}, the PV loop for a single cardiac cycle of the LV model is plotted for several timestep sizes. As the timestep size is decreased from $10^{-3} \unit{s}$ to $10^{-6} \unit{s}$, the PV loop converges. As shown in the figure, the PV loop for $10^{-3} \unit{s}$ is very close to the PV loop for $10^{-6} \unit{s}$, except for the maximum pressure, for which the difference is only about 1\%.

\subsection{Computational cost} \label{sect:cost}

\begin{figure}
    \centering
    \includegraphics[width=0.9\columnwidth]{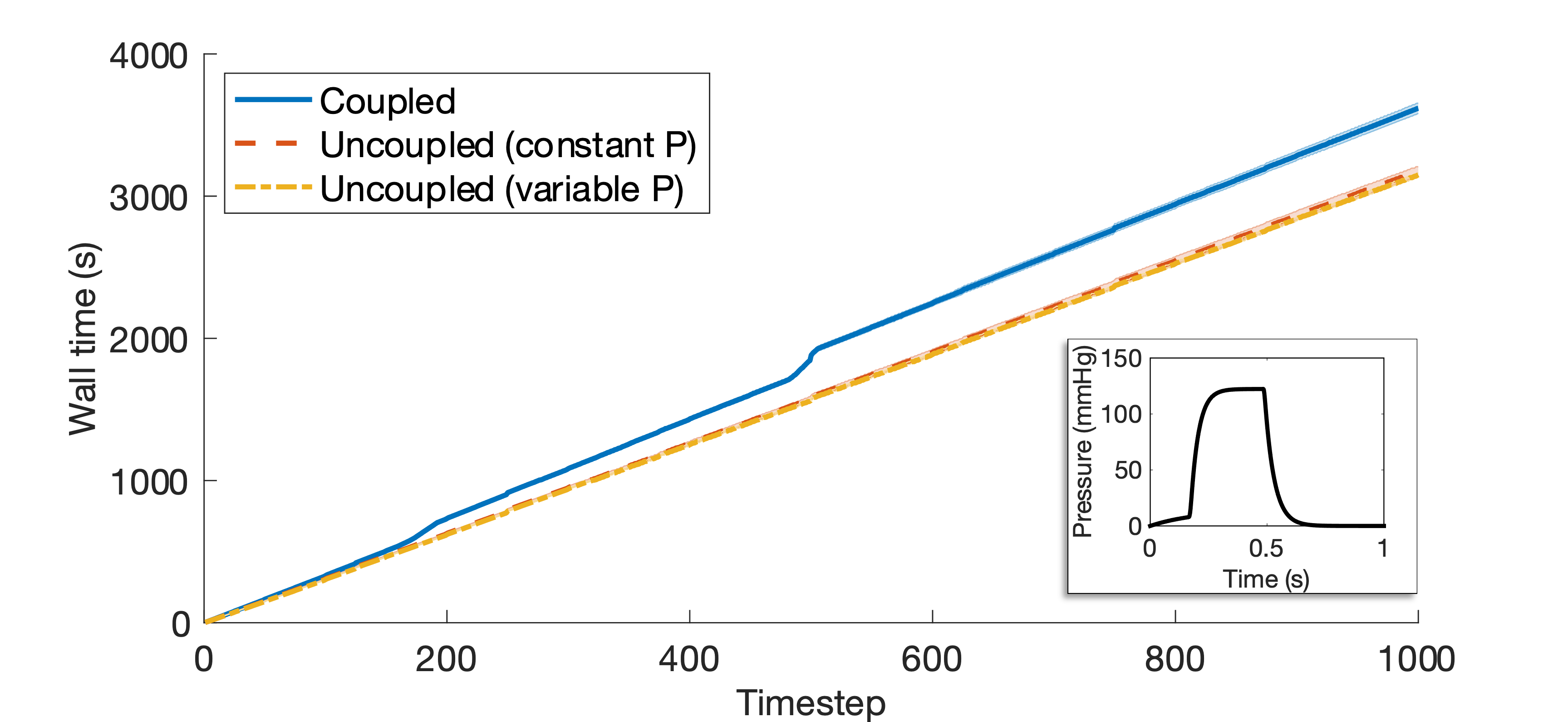}
    \caption{Comparison of simulation time for 0D-coupled versus uncoupled simulations. ``Coupled" is identical to that in Fig. \ref{fig:LV_sim_setup}, except it is run for only 1 cardiac cycle (1000 timesteps).
    ``Uncoupled (constant P)" is the same except the endocardial pressure is a constant 1500 \unit{Pa} (11.25 \unit{mmHg}). ``Uncoupled (variable P)" is the same except the endocardial pressure is given by a prescribed time-varying pressure curve, shown in the inset plot. All simulations were run in parallel with 4 processors. Each simulation was run 5 times, and the mean $\pm$ standard deviation of those samples are plotted. 
}
\label{fig:time_test}
\end{figure}
In this section, wall time is used to compare the computational cost of the 0D-coupled idealized LV simulation (Section \ref{sect:LV_sim}) with two similar but uncoupled simulations -- the first with a constant endocardial pressure, and the second with a time-varying endocardial pressure. All simulations were run for one cardiac cycle (1000 timesteps) with a timestep size of $10^{-3}\unit{s}$. All simulations were run in parallel with 4 CPUs of an Intel Gold 5118 2.3 GHz processor, and each case was run 5 times. The mean wall time $\pm$ standard deviation was computed for each set of simulations. As seen in Fig. \ref{fig:time_test}, after one cardiac cycle, the coupled simulation is approximately 15\% slower than the uncoupled simulations. This is due to the extra 0D solver communication and computation that is performed at each Newton iteration of the 3D solver. Note also that around timesteps 200 and 500, the coupled simulation slows down. These timesteps correspond to the isovolumic phases when both valves are closed, where our method experiences somewhat greater difficulty in convergence.

\subsection{Inflation of a spherical shell through its limit point} \label{sect:sphere_inflation}

\begin{figure}
    \centering
    \includegraphics[width=\columnwidth]{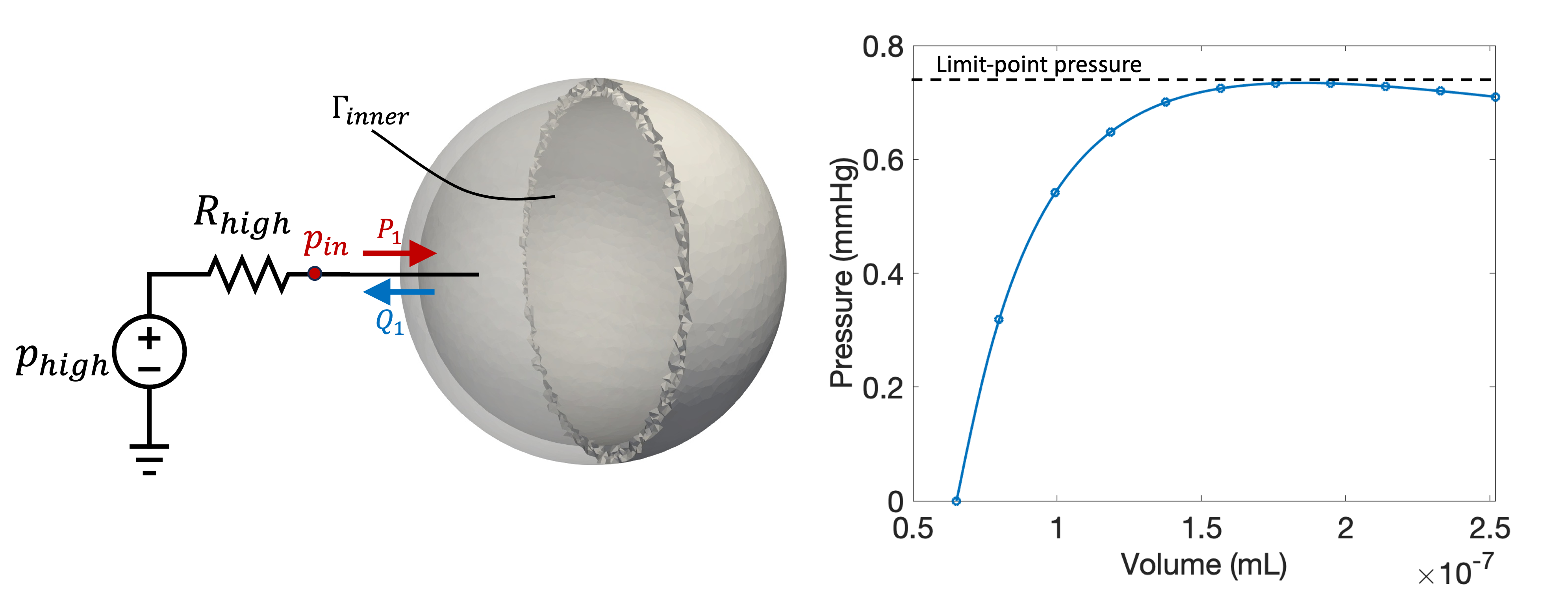}
    \caption{Left: A 3D thick spherical shell is coupled to a constant current 0D fluid. The inner surface of the shell is the coupled Neumann boundary. The 0D fluid model has a high-pressure source and a high resistance, which produces an approximately constant flow rate. Right: The pressure-volume relation for the spherical shell is plotted. Markers are placed every 100 ms. Inflating the sphere at roughly a constant rate of change of volume allows the limit point (where $dP/dV = 0$) to be traversed.
}
    \label{fig:sphere_inflation}
\end{figure}

An interesting feature of the coupling framework is that it also permits the investigation of so-called limit point problems in structural mechanics.
In this example, a simple 0D LPN is used to inflate a thick-walled sphere at approximately a constant flow rate (Fig. \ref{fig:sphere_inflation}). The coupled problem for this example may be stated in the structure of Eq. \eqref{eq:abstract_coupled_problem} as
\begin{equation}
    \begin{cases}
    p_{in} = p_{high} + Q_1 R_{high},
    \\
    \\
    \mathcal{P}^{3D,struct}(\mathbf{u}, \mathbf{x}, t) = \mathbf{0},
    \\
    \text{Initialize with zero displacement and velocity},
    \\
    \boldsymbol{\sigma} \cdot \mathbf{n} = -P_1 \mathbf{n}, \text{ on } \Gamma_{inner},
    \\
    \\
   P_1(t) = w_1(t) = p_{in}(t),
    \\
    Q_1(t) = \int_{\Gamma_{inner}(t)} \mathbf{u}(t) \cdot \mathbf{n}(t) d\Gamma.
    
    \end{cases}
\end{equation}
The LPN is the simplest ``constant current" circuit, consisting of a large resistance $R_{high}$ and a large pressure source $p_{high}$, chosen together to yield an approximately constant flow rate $-Q_1$ into the sphere. In this case, the LPN is described by a single algebraic equation for the pressure inside the sphere, $\mathbf{w} = [p_{in}]^T$. The coupled 0D Dirichlet boundary data are $\mathbf{Q} = [Q_1]^T$. The uncoupled 0D Neumann boundary data are $\mathbf{p} = [p_{high}]^T$, and there are no uncoupled 0D Dirichlet boundary data $\mathbf{q}_u = []$. By the coupling conditions, the pressure $P_1$ applied to the inner surface of the sphere, $\Gamma_{inner}$, is equal to $p_{in}$.  

The sphere has an inner radius of $25 \unit{\micro\meter}$ and a thickness of $2.5 \unit{\micro\meter}$, and it is composed of a neo-Hookean material with material constant $C_1 = 3 \unit{\kPa}$. With these parameters, the spherical shell displays limit point behavior. Specifically, as the sphere is inflated in a quasi-static process, the pressure increases, reaches a maximum (the limit point, where $dP/dV = 0$), then decreases \cite{anssari2022modelling, beatty1987topics} (Fig. \ref{fig:sphere_inflation} right). Similar concave down diastolic pressure-volume relations have been observed in embryonic chick hearts \cite{taber1992cardiac} and embryonic zebrafish hearts \cite{salehin2021assessing}, and limit point behavior has been proposed as a possible explanation. This example represents an idealized model of this phenomena in small-scale, embryonic hearts.

Limit point behavior presents a challenge to standard simulation methods. To simulate inflation of the sphere, one would typically apply an increasing pressure on the inner surface at constant increments. Unfortunately, this approach would never be able to capture the descending portion of the PV curve; once the applied pressure exceeds the limit point pressure, Newton's method will diverge and the simulation will crash because there is no static equilibrium configuration at that pressure. Techniques exist to traverse the limit point (and similar phenomena like snap-through behavior), the most common being the arc-length method \cite{crisfield1983arc}.

As shown in Fig. \ref{fig:sphere_inflation}, the present coupling framework allows us to traverse this limit point. There is in fact a connection between this coupling method and the arc-length method. In both cases, the load magnitude (pressure) is treated as unknown, and the equations corresponding to the mechanics problem are augmented by additional equations required to determine the load increment or load scaling factor. Actually, ANM was originally applied to solve coupled nonlinear systems arising from the arc-length method \cite{chan1985approximate}. We note that a standard monolithic 3D-0D coupling approach, such as in \cite{hirschvogel2017monolithic}, should also be able to traverse limit points in the same manner, but to the best of our knowledge, no previous studies have applied it to limit point problems.

\subsection{3D fluid - 0D fluid example} \label{sect:3dfluid_sim}
\begin{figure}
    \centering
    \includegraphics[width=\textwidth]{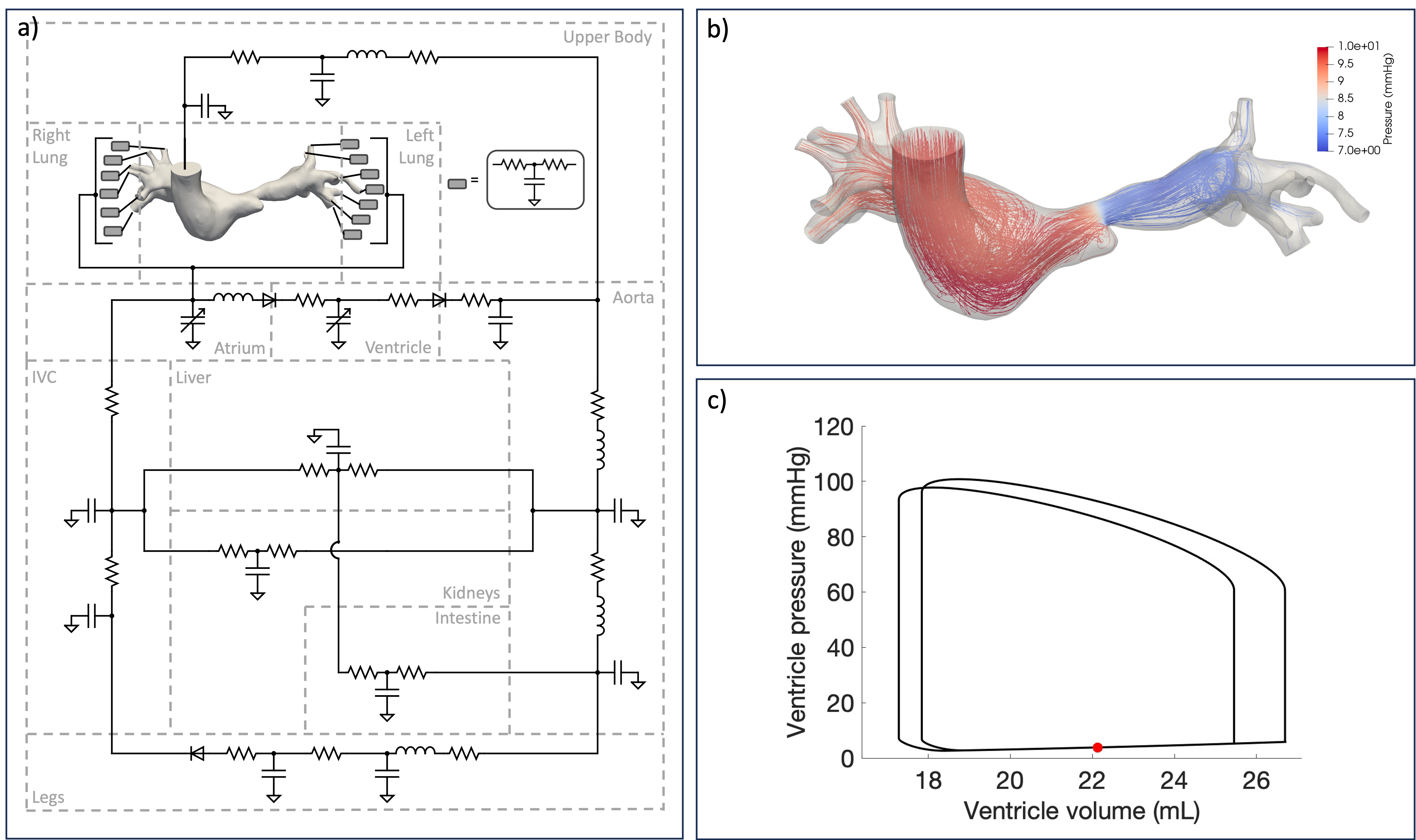}
    \caption{a) A 3D model of the pulmonary arteries is coupled to a closed-loop LPN model of the circulatory system \cite{schiavazzi2015hemodynamic}.
    b) Streamlines of blood flow, colored by pressure, in the pulmonary arterial model during mid-diastole (see red dot in PV loop panel c). c) The PV loop over two cardiac cycles of the LPN ventricle, which is modeled as a time-varying capacitor.}
    \label{fig:glenn_sim}
\end{figure}

For completeness, and to illustrate the generality of the coupling scheme to both 3D structure - 0D fluid problems and 3D fluid - 0D fluid problems, we reproduce a simulation similar to those in \cite{schiavazzi2015hemodynamic} using the current coupling scheme. In this work, the authors applied the original 3D fluid - 0D fluid coupling of \cite{moghadam2013modular} to simulate blood flow in the pulmonary arteries coupled to a complex, closed-loop model of the cardiovascular system. They used this model to investigate the hemodynamic effects of left pulmonary artery stenosis after the stage II superior cavo-pulmonary connection (SCPC) surgery. The model is illustrated in Fig. \ref{fig:glenn_sim}. The parameters of the model are taken from Table E1 P2 in \cite{schiavazzi2015hemodynamic}. 

The coupled problem for this example may be stated as,

\begin{equation}
    \begin{cases}
    \mathcal{P}^{0D}(\mathbf{w}, t, [\mathbf{q}_{u}, \mathbf{Q}], \mathbf{p}) = \mathbf{0} \text{ (equations not provided)},
    \\
    \text{Initial conditions},
    \\
    \\
    \mathcal{P}^{3D,fluid}(\mathbf{u}, p, \mathbf{x}, t) = \mathbf{0},
    \\
    \text{Initial conditions},
    \\
    \text{(Uncoupled) no-slip boundary conditions on walls,}
    \\
    \boldsymbol{\sigma} \cdot \mathbf{n} = -P_1 \mathbf{n}, \text{ on } \Gamma_{h_c}^{(1)},
    \\
    \hspace{0.75in} \vdots
    \\
    \boldsymbol{\sigma} \cdot \mathbf{n} = -P_{14} \mathbf{n}, \text{ on } \Gamma_{h_c}^{(14)},
    \\
    \\
   P_1(t) = (\mathbb{P}\mathbf{w}(t))_1,
   \\
    \hspace{0.75in} \vdots
    \\
    P_{14}(t) = (\mathbb{P}\mathbf{w}(t))_{14},
    \\
    Q_1(t) = \int_{\Gamma_{h_c}^{(1)}(t)} \mathbf{u}(t) \cdot \mathbf{n}(t) d\Gamma,
    \\
    \hspace{0.75in} \vdots
    \\
    Q_{14}(t) = \int_{\Gamma_{h_c}^{(14)}(t)} \mathbf{u}(t) \cdot \mathbf{n}(t) d\Gamma,
    
    \end{cases}
\end{equation}
where $\Gamma_{h_c}^{(1)}, \hdots, \Gamma_{h_c}^{(14)}$ are the 14 faces (1 inlet + 13 outlet) of the pulmonary model, which are associated with 14 coupling pressures $P_1 \hdots P_{14}$ and flow rates $Q_1 \hdots Q_{14}$.

For the sake of clarity, when applied to 3D fluid - 0D fluid problems, the coupling code used in the present work is functionally identical to the code used in \cite{schiavazzi2015hemodynamic}. We include this section to emphasize that a 3D fluid - 0D fluid problem can be solved using the current unified coupling scheme; however, the method when applied to the 3D fluid - 0D fluid problem and the results in this section are not novel. 

Fig. \ref{fig:glenn_sim} shows the 3D fluid and 0D fluid models, as well as the simulation results, which include a visualization of blood flow in the pulmonary arteries and the PV loop of the LPN ventricle component. 
Other examples of 3D fluid - 0D fluid coupling using the present coupling method can be found in \cite{moghadam2013modular, yang2019evolution, schwarz2021hemodynamic}.

\section{Discussion} \label{sect:discussion}
We present a general framework for coupling 3D models of cardiovascular components -- cardiac structures like the LV, as well as vascular blood flow models like the pulmonary arteries -- to 0D models of the circulatory system. The novel aspects of this work are the extension of the modular framework of \cite{moghadam2013modular} to the 3D structure - 0D fluid problem, a demonstration of the scheme's effectiveness in 0D-coupled cardiac mechanics applications, especially in its ability to avoid the balloon dilemma and capture isovolumic cardiac phases with ease, and finally a new derivation of the coupling scheme based on ANM, which provides important mathematical underpinnings and clearly shows the close connection to the monolithic coupling approach. 

Previous approaches to the coupling problem use either a partitioned or monolithic strategy, and the advantage of the present coupling is that it adopts a hybrid strategy, retaining the attractive features of both. In a monolithic approach, such as \cite{hirschvogel2017monolithic, pfaller2019importance}, modifying the 0D model can be cumbersome, requiring a detailed understanding of the 3D solver structure. In particular, one must typically derive new expressions for the off-diagonal tangent blocks in the monolithic (Newton) linear system Eq. \eqref{eq:monolithic_system}. In our approach, one can modify the 0D model without touching the 3D solver at all, and vice versa. As shown in Section \ref{sect:methods}, this is achieved by applying Schur complement reduction to Eq. \eqref{eq:monolithic_system}, which separates the 3D and 0D computations, then using the 0D fixed point operator to approximately solve the modified linear system.

The present approach corresponds to the Neumann coupling in \cite{moghadam2013modular}, in which a Neumann (pressure) boundary condition is imposed on the 3D domain, while a Dirichlet boundary condition is imposed on the 0D domain. In \cite{moghadam2013modular}, the authors also described a Dirichlet coupling, in which a Dirichlet boundary condition is imposed on the 3D domain, while a Neumann boundary condition is imposed on the 0D domain. In this type of coupling, the shape of the velocity profile on the 3D Dirichlet surface must be chosen. While applicable to fluid simulations, in which one can reasonably assume a parabolic or flat profile, this is an inconvenient limitation for structural mechanics simulations, since in most cases, the shape of the velocity profile cannot reasonably be assumed \textit{a priori}. In this paper, we only consider Neumann coupling. 

A different kind of Dirichlet coupling for the 3D structure - 0D fluid problem may be achieved by imposing a volume for the 3D structure, enforced by augmenting the structural equations with a volume constraint, as in recent works \cite{regazzoni2022cardiac, regazzoni2022stabilization, augustin2021computationally}. The augmented equations are then solved using Newton-like methods. Our Neumann coupling approach does not require such a volume constraint. Furthermore, we do not require a stabilization term like that introduced in the loosely coupled approach of \cite{regazzoni2022stabilization}, which was needed to eliminate unphysical oscillations near the isovolumic phases. Our approach automatically captures the isovolumic phases without unphysical oscillations, and without any adaptive/refined time stepping or any other special numerical techniques.

The authors of \cite{moghadam2013modular} also identified three variations of the method. In the ``implicit" method, the $\mathbf{M}$ matrix is updated each Newton iteration. In the ``semi-implicit" method, $\mathbf{M}$ is only computed once at the start of the simulation. In the ``explicit" method, the entire 0D contribution to the 3D tangent matrix is ignored. In \cite{moghadam2013modular}, the three variations were compared, and the semi-implicit method was found to provide the best balance of stability and cost-effectiveness. In this work, we described and used only the implicit method. Both the semi-implicit and explicit methods were unstable for our coupled LV simulation, likely due to the significant magnitude and drastic change in resistance produced by the 0D valves. We leave a detailed comparison of the three variations for coupled cardiac mechanics simulations to future work.

In multiple examples, we have shown the effectiveness of the present coupling method. To illustrate the application to cardiac mechanics modeling, we simulated an idealized LV coupled to an open-loop LPN over 10 cardiac cycles (Section \ref{sect:LV_sim}). The resulting PV loop spans physiological ranges of pressure and volume, shows a clear distinction of the four cardiac phases, and reaches a limit cycle. Critically, the isovolumic phases of the cardiac cycle are automatically captured without numerical stability issues. These results provide significant verification of the effectiveness of the coupling framework for cardiac mechanics simulations. Using the coupled LV model, we made preliminary assessments of the temporal convergence and cost of our method. In Section \ref{sect:temporal_convergence}, we showed that the coupling algorithm is stable and converges for small timesteps.  
In Section \ref{sect:cost}, we showed that coupled simulations are moderately more expensive than uncoupled simulations (15\% longer for one cardiac cycle). The additional cost comes from communication with the 0D solver and a few assembly operations each Newton iteration of the 3D solver. This relatively small increase in cost is well-warranted, given the increased physiological fidelity afforded by the coupling. 
In the spherical shell inflation case (Section \ref{sect:sphere_inflation}), which represents a simplified model of an embryonic heart, we showed 3D-0D coupling can be used to traverse a limit point and discussed its relationship to the arc-length method. Finally, we simulated blood flow in a pulmonary arterial model coupled to a closed-loop LPN of the circulatory system, demonstrating that our coupling also applies to the 3D fluid - 0D fluid problem (Section \ref{sect:3dfluid_sim}).

\paragraph{Limitations and Future Work}

With respect to the present coupling algorithm, we identify important areas for future investigation and development. 
We will first analyze the temporal convergence of the method, including the order of convergence as well as the effect of the 3D and 0D time-stepping schemes on overall convergence. In addition, under certain assumptions met in our case, ANM should converge quadratically \cite{chan1985approximate}. We intend to verify this quantitatively.

We are also interested in applying the stabilized structural mechanics formulation presented in \cite{liu2018unified} for cardiac mechanics simulations. 
This formulation treats both fluid and structural problems under a unified continuum modeling framework in which pressure is a primitive variable. It is highly effective for incompressible solids, which cardiac tissues are often assumed to be. Generally, we plan on applying the coupling in more complex cardiac mechanics simulations, for example with advanced myocardial constitutive models \cite{zhang2023simulating} or using 4-chamber anatomies \cite{jafari2019framework}.

This framework can be extended to couple a 3D model to a 1D model of blood circulation \cite{blanco2007unified}, which, unlike a 0D model, can recapitulate wave propagation phenomenon in arteries \cite{boileau2015benchmark}. On the 3D side, the expressions will be identical. One only needs to define the effect of a flow rate boundary condition on the 1D system, and define how to extract pressure from the 1D system to communicate back to the 3D system.

Finally, more recent coupling ideas similar to ANM have been summarized and proposed in \cite{yeckel2009approximate}, in which each solver is treated as a blackbox defined only through its fixed point iteration operator. These algorithms are even more modular than the present coupling scheme, while retaining the quadratic convergence properties of Newton's method. Future work may aim to apply these ideas to the 3D-0D coupling discussed here.

\section{Conclusion}\label{sect:conclusion}
In this work, a unified and modular framework for 3D-0D coupling in cardiovascular simulations is introduced. The algorithm, originally described in \cite{moghadam2013modular} for the 3D fluid - 0D fluid problem, is extended to solve the closely-related 3D structure - 0D fluid problem, showing that both problems can be treated uniformly within the same mathematical formulation. Through multiple examples, the effectiveness of the coupling algorithm is demonstrated. Notably, we construct a 0D-coupled idealized LV model that produces a physiological pressure-volume loop and effectively captures the isovolumic cardiac phases without additional numerical treatment. We also provide a new derivation using ANM, which reveals the present coupling scheme's connection to the monolithic Newton coupling approach. This hybrid coupling strategy combines the stability of monolithic approaches with the modularity and flexibility of partitioned approaches, with relatively small additional computational cost compared to uncoupled simulations. Overall, this work provides a robust, flexible, and efficient method for modeling the circulatory system in cardiovascular simulations of tissue mechanics and blood flow.

\section*{Acknowledgements}
Funding for this research was provided by the National Institutes of Health grants 5R01HL159970-02 and 5R01HL129727-06. The authors wish to thank Dr. Erica Schwarz, Dr. Fannie Gerosa, and Reed Brown for their helpful discussions during this project and for their comments on this paper.

\begin{appendices}
\section{4th Order Runge-Kutta scheme}
\label{sect:RK4_scheme}
RK4 applied to the ODE system Eq. \eqref{eq:0D_ode_with_q} reads
\begin{equation}  \label{eq:rk4}
\begin{aligned}
    &\mathbf{k}_1 = \mathbf{f}\Big(\mathbf{y}_n,
    \mathbf{z}_n, t_n, \mathbf{q}(t_n), \mathbf{p}(t_n))
    \Big)
    ,
    \\
    &\mathbf{k}_2 = \mathbf{f}\Big(
  \mathbf{y}_n + \mathbf{k}_1 \frac{\Delta t}{3}, \mathbf{z}_n, t_n + \frac{\Delta t}{3}, \mathbf{q}(t_n+\frac{\Delta t}{3}), \mathbf{p}(t_n+\frac{\Delta t}{3})
  \Big)
    ,
    \\
    &\mathbf{k}_3 = \mathbf{f}\Big(
  \mathbf{y}_n - \mathbf{k}_1 \frac{\Delta t}{3} + \mathbf{k}_2 \Delta t, \mathbf{z}_n, t_n + \frac{2 \Delta t}{3}, \mathbf{q}(t_n+\frac{2\Delta t}{3}), \mathbf{p}(t_n+\frac{2\Delta t}{3})
  \Big)
  ,
  \\
  & \mathbf{k}_4 = \mathbf{f} \Big(
  \mathbf{y}_n + \mathbf{k}_1 \Delta t - \mathbf{k}_2 \Delta t + \mathbf{k}_3 \Delta t, \mathbf{z}_n, t_n + \Delta t, \mathbf{q}(t_n + \Delta t),
  \mathbf{p}(t_n + \Delta t)
  \Big)
  ,
  \\
  & \mathbf{y}_{n+1} = \mathbf{y}_{n} + \frac{\mathbf{k}_1 +3\mathbf{k}_2 +3\mathbf{k}_3 +\mathbf{k}_4}{8} \Delta t
  .
\end{aligned}
\end{equation}
The algebraic variables are then determined by solving Eq. \eqref{eq:0D_algebraic_with_q} for $\mathbf{z}_{n+1}$ with the updated differential variables $\mathbf{y}_{n+1}$,
\begin{equation*} 
     \mathbf{g}(\mathbf{y}_{n+1}, \mathbf{z}_{n+1}, t_n + \Delta t, \mathbf{q}(t_n + \Delta t), \mathbf{p}(t_n + \Delta t)) = \mathbf{0}.
\end{equation*}
Rearranging this equation for $\mathbf{z}_{n+1}$ yields for some function $\tilde{\mathbf{g}}(\mathbf{y}, t, \mathbf{q}, \mathbf{p})$,
\begin{equation} \label{eq:z_np1_eqn}
    \mathbf{z}_{n+1}
    = 
    \tilde{\mathbf{g}}(\mathbf{y}_{n+1}, t_n + \Delta t, \mathbf{q}(t_n + \Delta t), \mathbf{p}(t_n + \Delta t))
    .
\end{equation}

Note that we assume the flow and pressure boundary forcings, $\mathbf{q}$ and $\mathbf{p}$, are known functions of time. In the coupled problem, $\mathbf{q}$ is composed of a prescribed uncoupled component  $\mathbf{q}_u$ and a coupled component $\mathbf{Q}$,
\begin{equation*}
    \mathbf{q}(t) = [\mathbf{q}_u(t), \mathbf{Q}(t)]
\end{equation*}
Since $\mathbf{q}_u$ is prescribed, its value is known at any time $t$. 
The value of $\mathbf{Q}$, on the other hand, is obtained from the 3D system, and its variation with time may be approximated by interpolating between its values at timestep $n$ and $n+1$,
\begin{equation*}
    \mathbf{Q}(t_n + h) = \mathbf{Q}_n + (\mathbf{Q}_{n+1} - \mathbf{Q}_n) \frac{h}{\Delta t},
\end{equation*}
where $\mathbf{Q}_n$ and $\mathbf{Q}_{n+1}$ are calculated from the 3D degrees of freedom $\boldsymbol{\Phi}_n$ and $\boldsymbol{\Phi}_{n+1}$, respectively.

From Eqs. \eqref{eq:rk4} and \eqref{eq:z_np1_eqn}, we can identify the fixed point operator corresponding to this 0D time integration scheme as
\begin{multline*}
    F^{0D, RK4}(\mathbf{w}_{n+1}, \boldsymbol{\Phi}_{n+1}; \mathbf{w}_{n}, \boldsymbol{\Phi}_{n}) 
    =
    \begin{bmatrix}
        \mathbf{y}_{n+1}
        \\
        \mathbf{z}_{n+1}
    \end{bmatrix}
    \\
    = 
    \begin{bmatrix}
        \mathbf{w}_n 
        + 
        \displaystyle \frac{\mathbf{k}_1 +3\mathbf{k}_2 +3\mathbf{k}_3 +\mathbf{k}_4}{8} \Delta t
        \\
        \tilde{\mathbf{g}}
        \Big(
        \mathbf{w}_n 
        + 
        \displaystyle \frac{\mathbf{k}_1 +3\mathbf{k}_2 +3\mathbf{k}_3 +\mathbf{k}_4}{8} \Delta t
        , 
        t_n + \Delta t, \mathbf{q}(t_n + \Delta t), \mathbf{p}(t_n + \Delta t)
        \Big)
    \end{bmatrix}
    .
\end{multline*}
\section{Calculation of flow rate} \label{sect:flow_rate}
For a coupled Neumann boundary condition, the 3D model must compute a flow rate $Q$ and send this value to the 0D domain. For 3D fluid - 0D fluid coupling (e.g., blood flow in the aorta coupled to Windkessel LPN), this is easily computed by integrating the normal component of the velocity over the coupled surface, $Q = \int_{\Gamma} \mathbf{u} \cdot \mathbf{n} d\Gamma$. For structural simulations, we must specify the context. In cardiac mechanics, we are typically interested in modeling a chamber or chambers of the heart. In this context, we are interested in the flow rate of blood into or out of a chamber, which is identical to the rate of change of volume of the chamber. We use this interpretation to compute flow rate in the 3D structure - 0D fluid problem.

Let $\Gamma$ be a closed surface, for example a sphere. The rate of change of volume enclosed by $\Gamma$ can be computed using the Reynolds Transport Theorem, which states that for an arbitrary scalar function $f$ defined over a time-varying region $\Omega(t)$ with boundary $\Gamma(t)$
\begin{equation}
    \frac{d}{dt} \int_{\Omega(t)} f d\Omega = \int_{\Omega(t)} \frac{\partial f}{\partial t} d\Omega + \int_{\Gamma(t)} (\mathbf{u_b} \cdot \tilde{\mathbf{n}}) f d\Gamma,
\end{equation}
where $\mathbf{u_b}$ is the velocity of the boundary and $\tilde{\mathbf{n}}$ is the outward unit normal of the boundary. Taking $f = 1$ and observing that $\frac{\partial f}{\partial t} = 0$, we find
\begin{equation*}
    \frac{d}{dt} \int_{\Omega(t)} (1) d\Omega = \int_{\Omega(t)} 0 d\Omega + \int_{\Gamma(t)} (\mathbf{u_b} \cdot \tilde{\mathbf{n}}) (1) d\Gamma,
\end{equation*}
\begin{equation} 
    \frac{dV}{dt} = \int_{\Gamma(t)} (\mathbf{u_b} \cdot \tilde{\mathbf{n}}) d\Gamma.
\end{equation}
That is, the rate of change of volume enclosed by a surface $\Gamma$ is the velocity flux integral over that surface. In cardiac mechanics, it is more natural to use the inward surface normal $\mathbf{n} = -\tilde{\mathbf{n}}$ (pointing away from the cardiac tissue and into the blood pool on the endocardial surface, for example). In addition, a positive flow rate $Q$ out of the heart chamber is associated with a decrease in chamber volume. Thus, we may write
\begin{equation} \label{eq:Q_derivation}
    Q = -\frac{dV}{dt} = \int_{\Gamma(t)} (\mathbf{u_b} \cdot \mathbf{n}) d\Gamma.
\end{equation}
This is in fact the same integral as in blood flow simulations, except we must be careful to take the integral over the coupled surface in the current (deformed) configuration, $\Gamma(t)$.

It is important to note that Eq. \eqref{eq:Q_derivation} is valid only if $\Gamma$ is a closed surface. In general, this is not the case, for example in the left ventricle or biventricle models cut at the basal plane. For these models, in order to accurately calculate $\displaystyle \frac{dV}{dt}$, it is necessary to close $\Gamma$ with a ``cap" surface. Such a capping was done in this work and slightly modifies the expressions for the 0D residual and tangent contributions (see Section \ref{sect:capping}).

\section{Contribution of coupled surface to tangent matrix} \label{sect:0D_tangent_contribution_derivation}
Here we derive Eq. \eqref{eq:0D_tangent_contribution} from Eq. \eqref{eq:0D_tangent_contribution_abstract}. Eq. \eqref{eq:0D_tangent_contribution_abstract} reads
\begin{equation*}
    \begin{bmatrix}
        \mathbf{K}^{3D/0D}
    \end{bmatrix}
    _{n+1}^{(k)}
    =
    \begin{bmatrix}
    \displaystyle \frac{\partial \mathbf{R}^{3D}}{\partial \mathbf{P}_{n+1} } 
    \mathbb{P}
    \displaystyle \frac{\partial F^{0D}(\mathbf{w}_{n+1}, \mathbf{Q}_{n+1})}{\partial \mathbf{Q}_{n+1}}
    \displaystyle \frac{\partial \mathbf{Q}_{n+1}}{\partial \mathbf{U}_{n+1}}
    \displaystyle \frac{\partial \mathbf{U}_{n+1}}{\partial \boldsymbol{\Phi}_{n+1}}
    \end{bmatrix}
    _{n+1}^{(k)}
    .
\end{equation*}
We consider each term individually. 

\begin{itemize}
    \item 
    The residual $\mathbf{R}^{3D}$ depends on the coupling pressures $\mathbf{P}_{n+1}$ through an equation like Eq. \eqref{eq:0D_residual_contribution},
    \begin{equation*}
   R_{Ai}^{3D} = \text{other terms} + \sum_{j=1}^{n^{cBC}}\int_{\Gamma_{h_c}^{(j)}} N_A P_{n+1,j} n_i d\Gamma
   .
    \end{equation*}
    Thus,
    \begin{equation}
        \frac{\partial (R_{Ai}^{3D})}{\partial P_{n+1,j}} = \int_{\Gamma_{h_c}^{(j)}} N_A  n_i d\Gamma.
    \end{equation}

    Note that if the 3D residual is composed of momentum and continuity components (e.g., for Navier-Stokes), then this term is multiplied by $\mathbb{I}_m = \begin{bmatrix}
        \mathbf{1} \\
        \mathbf{0}
    \end{bmatrix}$, which is simply a vector with 1s corresponding to momentum rows and 0s corresponding to continuity rows. This places the tangent contribution in the momentum block row.

    \item In Section \ref{sect:0D_tangent_contribution}, we described how $\mathbb{P} \displaystyle \frac{\partial F^{0D}(\mathbf{w}_{n+1}, \mathbf{Q}_{n+1})}{\partial \mathbf{Q}_{n+1}}$ is computed in a finite difference manner (Eq. \eqref{eq:M_matrix}). This term is denoted by the resistance-like matrix $M_{ij}$.

    \item From Eq. \eqref{eq:Q3D}, we have
    \begin{equation*}
        Q_{n+1,i} = \sum_A \int _{\Gamma_{h_c}^{(i)}} N_A U_{n+1,Ak} n_k d\Gamma.
    \end{equation*}
    Thus,
    \begin{equation}
        \frac{\partial Q_{n+1,i}}{\partial U_{n+1,Ak}} = \int_{\Gamma_{h_c}^{(i)}} N_A n_k d\Gamma.
    \end{equation}

    \item Finally, we deal with the term $\displaystyle \frac{\partial \mathbf{U}_{n+1}}{\partial \boldsymbol{\Phi}_{n+1}}$. This term depends on the time discretization scheme. In our case, we use the generalized-$\alpha$ method and choose the nodal accelerations $\dot{\mathbf{U}}$ as our 3D unknowns (along with nodal pressures $\mathbf{\Pi}$ if the 3D is a fluid). The nodal velocities $\mathbf{U}$ are related to the nodal accelerations $\dot{\mathbf{U}}$ by Eq. \eqref{eq:genalpha_U}
    \begin{equation*}
        \mathbf{U}_{n+1} = \mathbf{U}_n + \Delta t \dot{\mathbf{U}}_{n} + \gamma \Delta t (\dot{\mathbf{U}}_{n+1} - \dot{\mathbf{U}}_n).
    \end{equation*}
    Thus,
    \begin{equation}
        \frac{\partial U_{n+1,Ai}}{\partial \Phi_{n+1,Bj}} = \frac{\partial U_{n+1,Ai}}{\partial \dot{U}_{n+1,Bj}} = \gamma \Delta t \delta_{AB} \delta_{ij}.
    \end{equation}

    Note that if the 3D degrees of freedom contain both acceleration and pressure components (e.g., for Navier-Stokes), then this term is multiplied by $\mathbb{I}_a = \begin{bmatrix}
        \mathbf{1} \\
        \mathbf{0}
    \end{bmatrix}$, which is simply a vector with 1s corresponding to acceleration rows and 0s corresponding to pressure rows. This places the tangent contribution in the acceleration block column.
\end{itemize}
Combining these terms yields the tangent matrix contribution Eq. \eqref{eq:0D_tangent_contribution}. In addition to this tangent matrix contribution, with a follower pressure load there should be additional tangent contribution due to the fact that $\Gamma_{h_c}^{(i)}$ and $n_i$ change with the deformation. However, this term is no different than the required term for an uncoupled surface and is unrelated to our coupling method, so we do not include it here.

\section{Idealized LV simulation parameters} \label{sect:sim_params}
We list the parameter values for the idealized LV coupled to open-loop LPN simulation, shown in Fig. \ref{fig:LV_sim_setup}. Most values are taken from \cite{pfaller2019importance}.

\begin{center}
\begin{longtable}{c  c  c  c}  \label{tab:LV_sim_params}
 Name & Parameter & Value & Unit \\ [0.5ex] 
 \hline\hline
\textbf{\textit{General mechanical}}  \\
 Tissue Density     & $\rho_0$      & $10^3$        & $\unit{kg/m^3}$ \\ 

 Viscosity          & $\eta$        & $100$        & $\unit{Pa \cdot s}$ \\

 Volumetric Penalty & $\kappa$      & $10^6$        & $\unit{Pa}$ \\

\\
 \textbf{\textit{Active stress}} \cite{pfaller2019importance} \\

 Contractility      & $\sigma_0$   &   $\num{8e4}$    & $\unit{Pa}$            \\

 Activation rate    & $\alpha_{max}$  &   $+5$      &  $\unit{1/s}$    \\
 
 Deactivation rate  & $\alpha_{min}$   &  $-30$      &  $\unit{1/s}$    \\
 
 Ventricular systole     & $t_{sys,v}$       &   $0.143$      &  $\unit{s}$    \\
 
 Ventricular diastole     & $t_{dias,v}$       &   $0.484$      &  $\unit{s}$     \\
 
 Steepness             & $\gamma$       &   $0.005$      &  $\unit{s}$    \\

\\
 \textbf{\textit{Passive myocardial tissue (HO model)}}                  \\
 
 Matrix             & $a$       &   $59.0$      &  $\unit{Pa}$    \\
 
                     & $b$       &   $8.023$      &  $\unit{-}$    \\
 
 Fiber             & $a_f$       &   $\num{18.472e3}$      &  $\unit{Pa}$     \\
 
                 & $b_f$       &   $16.026$      &  $\unit{-}$    \\
 
 Sheet             & $a_s$       &   $\num{2.481e3}$      &  $\unit{Pa}$    \\
  
                 & $b_s$       &   $11.12$       &  $\unit{-}$    \\
  
 Fiber sheet        & $a_{fs}$       &   $216$      &  $\unit{Pa}$     \\
  
                 & $b_{fs}$       &   $11.436$      &  $\unit{-}$    \\

\\
\textbf{\textit{Epicardial boundary condition}}                 \\

 Spring stiffness     & $k_{epi}$       &   $\num{1.0e8}$      &  $\unit{Pa/m}$    \\
  
 Dashpot viscosity     & $c_{epi}$       &   $\num{5.0e3}$      &  $\unit{Pa \cdot s/m}$    \\

\\
 \textbf{\textit{Basal boundary condition}}                \\
 
 Spring stiffness     & $k_{base}$       &   $\num{1.0e5}$      &  $\unit{Pa/m}$    \\
  
 Dashpot viscosity     & $c_{base}$       &   $\num{5.0e3}$      &  $\unit{Pa.s/m}$    \\

\\
\textbf{\textit{LPN}}     \\
 Proximal inertance     & $L_p$       &   $\num{1.3e5}$      &  $\unit{kg/m^4}$    \\
  
 Proximal capacitance    & $C_p$       &   $\num{7.7e-9}$      &  $\unit{m^4.s^2/kg}$   \\
  
 Distal capacitance      & $C_d$       &   $\num{8.7e-9}$      &  $\unit{m^4.s^2/kg}$    \\
  
 Proximal resistance     & $R_p$       &   $\num{7.3e6}$      &  $\unit{kg/m^4/s}$    \\
  
 Distal resistance        & $R_d$       &   $\num{1.0e8}$      &  $\unit{kg/m^4/s}$    \\
   
 Reference pressure        & $P_{ref}$   &   $0.0$      &  $\unit{Pa}$    \\
   
 Closed valve resistance    & $R_{max}$    &   $\num{1e9}$      &  $\unit{kg/m^4/s}$    \\
   
 Open valve resistance      & $R_{min}$       &   $\num{1e6}$     &  $\unit{kg/m^4/s}$    \\
   
 Valve steepness             & $k_p$       &   $\num{1e-3}$      &  $\unit{Pa}$    \\

 Baseline atrial pressure    & $p_{at,0}$       &   $6.0$      &  $\unit{mmHg}$    \\
   
 Atrial pressure amplitude   & $\Delta p_{at}$   &   $8.0$      &  $\unit{mmHg}$    \\
   
 Atrial systole             & $t_{sys,a}$       &   $0.0$      &  $\unit{s}$    \\
   
                            & $T_{sys,a}$       &   $0.2$      &  $\unit{s}$    \\

\\
\textbf{\textit{LPN initial conditions}} \\
 Initial ventricle pressure    & $p_v(0)$       &   $8.0$      &  $\unit{mmHg}$    \\

Initial proximal pressure  & $p_p(0)$       &   $61.8$      &  $\unit{mmHg}$    \\

Initial distal pressure   & $p_d(0)$       &   $59.7$      &  $\unit{mmHg}$    \\

Initial proximal flowrate   & $q_p(0)$       &   $38.3$      &  $\unit{mL/s}$    \\

\\
\textbf{\textit{Numerical integration}}     \\
Timestep size           &  $\Delta t$   &   $0.001$ &    $\unit{s}$      \\
Gen-$\alpha$ parameter  &  $\gamma$     & $0.5$     &    $\unit{-}$     \\
Gen-$\alpha$ parameter  &  $\alpha_f$     & $0.5$     &    $\unit{-}$   \\
Gen-$\alpha$ parameter  &  $\alpha_m$     & $0.5$     &    $\unit{-}$   \\
Gen-$\alpha$ parameter  &   $\beta$       & $0.25$    &     $\unit{-}$  \\
 \hline
\end{longtable}
\end{center}

\end{appendices}

\bibliographystyle{unsrtnat}
\bibliography{references}  






\end{document}